\documentclass[twocolumn,amsthm]{autart}

\usepackage{amsmath, amssymb, bbm, xspace}
\usepackage{amsfonts,mathrsfs}
\usepackage{mathtools}
\usepackage{color}
\usepackage{graphicx}
\usepackage{epsfig}
\usepackage{algorithm}
\usepackage{algorithmicx}
\usepackage[acronym]{glossaries}
\usepackage{subfigure}
\usepackage{cancel}
\usepackage{empheq}
\usepackage{placeins}
\usepackage{enumitem}
\usepackage{textgreek}
\usepackage{empheq}
\usepackage{url}
\usepackage{xcolor}
\usepackage{caption}
\usepackage{booktabs}

\newcommand{\bs}{\boldsymbol}
\newcommand{\mc}{\mathcal}

\renewcommand{\emph}{\textit}

\makeglossaries
\newacronym{GNEP}{GNEP}{generalized Nash equilibrium problem}
\newacronym{NE}{NE}{Nash equilibrium}
\newacronym{NEP}{NEP}{Nash equilibrium problem}
\newacronym{GNE}{GNE}{generalized Nash equilibrium}
\newacronym{v-GNE}{v-GNE}{variational \gls{GNE}}
\newacronym{ISS}{ISS}{Input-to-state-stable}
\newacronym{PPA}{PPA}{proximal-point algorithm}
\newacronym{PPPA}{PPPA}{preconditioned \gls{PPA}}
\newacronym{VI}{VI}{variational inequality}
\newacronym{GAE}{GAE}{generalized aggregative equilibrium}
\newacronym{v-GAE}{v-GAE}{variational \gls{GAE}}
\newacronym{KKT}{KKT}{Karush--Kuhn--Tucker}
\newacronym{FQNE}{FQNE}{firmly quasinonexpansive}
\newacronym{FNE}{FNE}{firmly nonexpansive}
\newacronym{FB}{FB}{forward-backward}

\renewcommand{\iff}{\; \Leftrightarrow \;}
\newcommand{\0}{\bs 0}
\def\1{{\bs 1}}

\def\argmin{\mathop{\rm argmin}}

\newcommand{\minimize}[1]{\displaystyle\minim_{#1}}
\newcommand{\minim}{\mathop{\hbox{\rm minimize}}}

\newcommand{\col}{\mathrm{col}}
\newcommand{\avg}{\operatorname{avg}}

\def\diag{\mathop{\hbox{\rm diag}}}

\def\Null{\mathop{\hbox{\rm null}}}

\def\Range{\mathop{\hbox{\rm range}}}
\def\range{\mathop{\hbox{\rm range}}}

\def\spose#1{\hbox to 0pt{#1\hss}}

\newcommand{\proj}{\mathrm{proj}}
\def\fix{\mathrm{fix}}
\def\gra{\operatorname{gra}}
\def\zer{\operatorname{zer}}
\def\dom{\operatorname{dom}}

\def\J{\mathrm{J}{}}
\newcommand{\Id}{\mathrm{Id}}
\def\nc{\mathrm{N}}
\def\H{\mc{H}}

\def\R{\mathbb{R}}
\def\N{\mathbb{N}}

\def\o{\omega}

\def\l{\lambda}

\def\Fc{\mathcal{F}}

\def\Ac{\mathcal{A}}

\def\Fc{\mathcal{F}}
\def\nbar{\skew2\bar n}

\def\xbs{{\bs x}}

\def\obs{{\bs \o}}
\def\lbs{{\bs \l}}

\def\vbs{{\bs v}}
\def\ubs{{\bs u}}

\def\obsbar{{\bar{\bs \o}}}

\def\Fa{\bs{F}_{\!  \textrm{a}}}
\def\Fi{\bs{\tilde{F}}_{\! i}}
\def\ha{{ \hphantom{{}={}} }}
\newcommand{\alphamax}{\alpha_{\textnormal{max}}}
\newcommand{\alphamaxtilde}{\tilde{\alpha}_{\textnormal{max}}}
\def\k{{k \in \N}}  

\DeclareSymbolFont{myletters}{OML}{ztmcm}{m}{it}
\DeclareMathSymbol{\uplambda}{\mathord}{myletters}{"15}

\def\QEDhereeqn{\eqno\let\eqno\relax\let\leqno\relax\let\veqno\relax\hbox{\QED}}
\def\QEDopenhereeqn{\eqno\let\eqno\relax\let\leqno\relax\let\veqno\relax\hbox{\QEDopen}}

\newlist{thmlist}{enumerate}{1}
\setlist[thmlist]{label=(\roman{thmlisti}), ref=\thethm(\roman{thmlisti}),noitemsep}
\newlist{lemlist}{enumerate}{1}
\setlist[lemlist]{label=(\roman{lemlisti}), ref=\thelem(\roman{lemlisti}),noitemsep}

\newtheorem{lem1}{Lemma}

\newtheorem{prop1}{Proposition}
\newtheorem{standing}{Standing Assumption}
\theoremstyle{definition}
\newtheorem{defn1}{Definition}
\newtheorem{exmp1}{Example}
\newtheorem{rem1}{Remark}
\newenvironment{proof1}[1][Proof]{\textbf{#1.} }

\begin{document}
	\begin{frontmatter}
		\title{Fast generalized Nash equilibrium seeking under partial-decision information\thanksref{footnoteinfo}}
		\thanks[footnoteinfo]{
			{This work is supported by the NWO under project OMEGA (613.001.702) and by the ERC under  project COSMOS (802348). Some material in this paper was partially presented at the 2020 IEEE Conference on Decision and Control.} }
		\author[aff1]{Mattia Bianchi},
		\author[aff2]{Giuseppe Belgioioso},
		\author[aff1]{Sergio Grammatico}

		\address[aff1]{Delft Center for Systems and Control, Delft University of Technology, The Netherlands}
		\address[aff2]{Automatic Control Laboratory, Swiss Federal Institute of Technology (ETH) Z\"{u}rich, Switzerland}


		\begin{abstract}
			We address the generalized Nash equilibrium seeking problem in a partial-decision information scenario, where each agent can only exchange information with some neighbors, although its cost function possibly depends on the strategies of all agents. The few existing methods  build on projected pseudo-gradient dynamics, and require either double-layer iterations or conservative conditions on the step sizes. To overcome both these flaws and improve efficiency, we design the first fully-distributed single-layer algorithms based on proximal best-response. Our schemes are fixed-step and allow for inexact updates, which is crucial for reducing the computational complexity.
			 Under standard assumptions on the game primitives, we establish convergence to a variational equilibrium (with linear rate for games without coupling constraints) by  recasting our algorithms as proximal-point methods, opportunely preconditioned to distribute the computation among the agents.
			 Since our analysis hinges on a restricted monotonicity property, we also provide new general results that significantly extend the domain of applicability of proximal-point methods.
			  Besides, our operator-theoretic approach  favors the implementation of provably correct acceleration schemes that can further improve the convergence speed.
			 Finally, the potential of our algorithms is demonstrated numerically, revealing much  faster convergence
			with respect to projected pseudo-gradient  methods   and validating our theoretical findings.
		\end{abstract}

	\end{frontmatter}

	\section{Introduction}
	Generalized games model the  interaction between self-interested decision makers, or agents, that aim at optimizing their individual, yet inter-dependent, objective functions, subject to shared constraints.
	This competitive scenario has received increasing attention with the spreading of
	networked systems,
	due to the numerous engineering applications, including  demand-side management in the smart grid \cite{Saad2012},  charging/discharging of electric vehicles \cite{Grammatico2017}, demand response in competitive markets \cite{Li_Chen_Dahleh_2015},
	and radio communication \cite{Palomar_Eldar_facchinei_pang_2009}.
	From a game-theoretic perspective, the challenge is   to assign the agents behavioral rules that  eventually ensure the attainment of a satisfactory equilibrium.

	A recent part of the literature focuses in fact on designing \emph{distributed} algorithms to seek a \gls{GNE}, a  decision set from which no agent has interest to unilaterally deviate \cite{FacchineiKanzow2010}, \cite{Yu_VanderSchaar_Sayed_2017}, \cite{BelgioiosoGrammatico_ECC_2018}, \cite{YiPavel2019},  \cite{ChenMingHongYi_Aut2021}. In these works, the computational effort is partitioned among the agents, but assuming that each of them has access to the decision of all the competitors (or to an aggregation value, in the case of aggregative games). Such an hypothesis, referred to as \emph{full-decision information}, generally requires the presence of a central coordinator that communicates with all the agents, which is impractical in some cases \cite{SwensonKarXavier_FictitiousPlay_2015}, \cite{FrihaufKrsticBasar_2012}. One example is the Nash--Cournot  competition model described in \cite{Koshal_Nedic_Shanbag_2016}, where the  profit of each of a group of firms depends not only on its own production, but also on the total supply, a quantity not directly accessible by any of the firms.  Instead, in this paper we  consider the so-called \emph{partial-decision information} scenario, where each agent
	 estimates the actions of all the competitors by relying only on the information exchanged with some neighbors over a communication network. Thus, the goal is to design \emph{fully-distributed} (namely,  center-free) algorithms, based  exclusively on peer-to-peer communication.

	 The partial-decision information setup has only been introduced very  recently. Most results consider non-generalized games (i.e., games without  shared constraints) \cite{Koshal_Nedic_Shanbag_2016}, \cite{Salehisadaghiani_Pavel_GOSSIP}, \cite{TatarenkoShiNedicTAC20}, \cite{SalehisadaghianiWeiPavel2019}. Even fewer algorithms can cope with the presence of coupling constraints  \cite{Pavel2018},  \cite{BelgioiosoNedicGrammatico2020}, \cite{GadjovPavel_aggregative_2019}, despite this extension arises naturally in most resource allocation problems \cite[§2]{FacchineiKanzow2010}, e.g., due to shared capacity limitations.
	 All the cited formulations resort to (projected) gradient and consensus-type dynamics, and are single-layer (i.e., they require a fixed finite number of communications per iteration). The main drawback is that, due to the partial-decision information assumption, theoretical guarantees are obtained only for small (or vanishing) step sizes, which significantly affect the speed of convergence. The only alternative available in literature  consists of double-layer algorithms, \cite{Lei_Shanbag_CDC2018}, \cite{PariseGentileLygeros_TCNS2020}, where the agents must communicate multiple (virtually infinite) times to reach consensus, before each update. An extensive communication requirement is however a performance bottleneck, as the communication time can overwhelm the time spent on local useful processing -- in fact,  this is  a common problem in parallel computing \cite{NIPS2019}.
	 Let alone the time lost in the transmission, sending large volumes of data on wireless networks results in a dramatically increased energetic cost.

	\emph{Contributions:} To improve speed and efficiency, we design the first fully-distributed single-layer \gls{GNE} seeking algorithms based on proximal best-response. For the sake of generality and mathematical elegance, we take here an operator-theoretic approach \cite{BelgioiosoGrammatico_L_CSS_2017}, \cite{YiPavel2019}, and  reformulate the \gls{GNE} problem  as that of finding a zero of a monotone operator. The advantage is that several fixed-point iterations are known to solve monotone inclusions  \cite[§26]{Bauschke2017}, thus providing a unifying framework to design algorithms and study their convergence. For instance,  the methods in \cite{Pavel2018},  \cite{BelgioiosoNedicGrammatico2020}, \cite{GadjovPavel_aggregative_2019}, were developed based on the (preconditioned) \gls{FB} splitting \cite[§26.5]{Bauschke2017}.  To enhance the convergence speed, we instead employ a \gls{PPA} \cite[Th.~28.1]{Bauschke2017}, which typically can tolerate much larger step sizes. Nonetheless, the design of distributed \gls{GNE} seeking \glspl{PPA} was elusive until now, because a direct implementation results in double-layer algorithms \cite{Scutari_ComplexGames_TIT2014}, \cite{YiPavel_PPP_2019}.
    The novelties of this work are summarized as follows:
	\begin{itemize}[nolistsep,topsep=0em,leftmargin=*]
		\item We propose the first  \gls{PPA} to compute a zero of a \emph{restricted} monotone operator, which significantly generalizes  classical results for maximally monotone operators. Differently from other recent extensions  \cite{ElFarouq_PseudomonotonePPP_2001}, \cite{Moudafi2020}, we also allow for set-valued resolvents and inexact updates, and we do not assume pseudomonotonicity or hypomonotonicity. This is a fundamental result of independent interest,  which we  exploit to prove convergence of our algorithms (§\ref{subsec:PPA});
		\item  We introduce a novel primal-dual proximal best-response  \gls{GNE}  seeking algorithm, which is the first non-gradient-based scheme for the
		partial-decision information setup.
		We derive our method as a \gls{PPA}, where we design a novel preconditioning matrix to distribute the computation and obtain a single-layer iteration.
		Under strong monotonicity and Lipschitz continuity of the game mapping, we prove global convergence with fixed step sizes, by exploiting restricted monotonicity properties.
		Convergence is retained even if the proximal best-response is computed inexactly (with summable errors), which is crucial for practical implementation.
		Differently from  \cite[Alg.~1]{Pavel2018},
		the step sizes
		can be chosen independently of a certain restricted strong monotonicity  constant.  In turn, not  only we allow for much larger steps, but parametric dependence is also improved: for instance, the bounds do not vanish when the number of agents grows, and the resulting convergence rate for non-generalized games is superior.
		 Moreover our scheme requires only one communication per iteration, instead of two
		 (§\ref{subsubsec:derivation}, {§\ref{subsec:rate}});
		\item We apply some acceleration schemes \cite{Iutzeler_Hendrickx2019} to our \gls{PPPA} and provide new  theoretical convergence guarantees.  We observe numerically that the  iterations needed to converge can be halved (§\ref{sec:acc});
		\item We tailor our method to efficiently solve aggregative games, by letting each agent keep and exchange an estimate of the aggregative value only, instead of an estimate of all the other agents' actions (§\ref{sec:aggregative});
		\item Via numerical simulations, we show that our \glspl{PPPA}
         significantly outperform the pseudo-gradient methods in  \cite{Pavel2018}, \cite{GadjovPavel2018} (the only other known fully-distributed, single-layer, fixed-step \gls{GNE} seeking schemes), not only in terms of number of iterations needed to converge (hence with a considerable reduction of the communication burden),  but also in terms of total computational cost (despite each agent must locally solve a strongly convex optimization problem, rather than a projection, at each step)  (§\ref{sec:numerics}).
	\end{itemize}

   Some preliminary results of this paper appeared in \cite{Bianchi_CDC20_PPP},  where we study only one special case for games without coupling constraints and with exact computation of the resolvent (and where we do not consider aggregative games or acceleration schemes), see §\ref{subsec:rate}.

	\emph{Basic notation}:  $\mathbb{N}$ is the set of natural numbers, including $0$.
	$\R$ ($\R_{\geq 0}$) is the set of (nonnegative) real numbers.
	$\0_q\in \R^q$ ($\1_q\in\R^q$) is a vector with all elements equal to $0$ ($1$); $I_q\in\R^{q\times q}$ is an identity matrix; the subscripts may be omitted when there is no ambiguity.
	For  a matrix $A \in \R^{p \times q}$, $[A]_{i,j}$ is the element on  row $i$ and column $j$; $\Null(A)\coloneqq \{x\in\R^q \mid Ax=\0_n\}$ and $\range(A)\coloneqq \{v\in\R^p \mid v=Ax, x\in\R^q \}$; $\| A \|_\infty$ is the maximum of the absolute row sums of $A$.
	If $ A=A^\top\in\R^{q\times q}$, $\uplambda_{\textnormal{min}}(A)=:\uplambda_1(A)\leq\dots\leq\uplambda_q(A)=:\uplambda_{\textnormal{max}}(A)$ denote its eigenvalues.
	$\diag(A_1,\dots,A_N)$ is the block diagonal matrix with $A_1,\dots,A_N$ on its diagonal. Given $N$ vectors $x_1, \ldots, x_N$,  $\col (x_1,\ldots,x_N ) \coloneqq  [ x_1^\top \ldots  x_N^\top ]^\top$. $\otimes$ denotes the Kronecker product.
	$\ell^1$  is the set of  absolutely summable sequences.

	\emph{Euclidean spaces}: Given a positive definite matrix $ \R^{q\times q} \ni P \succ 0$,  $\mc{H}_P\coloneqq (\R^q,\langle \cdot \mid \cdot \rangle _P)$ is the Euclidean space obtained by endowing $\R^q$ with the $P$-weighted inner product  $\langle x \mid y \rangle _P=x^\top P y$, and $\| \cdot \|_P$ is the associated norm; we omit the subscripts if $P=I$. Unless otherwise stated, we always assume to work in $\H=\H_I$.

	\emph{Operator-theoretic background}:
	A set-valued  operator $\mc{F}:\R^q\rightrightarrows \R^q$ is characterized by its graph
	$\gra (\mc{F})\coloneqq \{(x,u) \mid u\in \mc{F}(x)\}$. $\dom(\mc{F})\coloneqq \{x\in\R^q| \mc{F}(x)\neq \varnothing \}$,
	$\fix\left( \mathcal{F}\right) \coloneqq  \left\{ x \in \R^q \mid x \in \mathcal{F}(x) \right\}$ and $\zer\left( \mathcal{F}\right) \coloneqq  \left\{ x \in \R^q \mid 0 \in \mathcal{F}(x) \right\}$ are the domain, set of fixed points and set of zeros, respectively. $\mc{F}^{-1} $ denotes the inverse operator of $\mc{F}$, defined as $\gra (\mc{F}^{-1})=\{(u,x)\mid (x,u)\in \gra(\mc{F})\}$.
	$\mathcal{F}$ is
	($\mu$-strongly) monotone in $\mc{H}_P$ if $\langle u-v \mid x-y\rangle_P \geq 0$ ($\geq \mu\|x-y\|^{2}_P$) for all $(x,u)$,$(y,v)\in\gra(\Fc)$;
	 we omit the indication ``in $\mc{H}_P$'' whenever $P=I$.
	$\Id$ is the identity operator. ${\rm J}_{\mathcal{F} }\coloneqq (\Id + \mathcal{F} )^{-1}$ denotes the resolvent operator of $\mathcal{F} $. For a function $\psi: \R^q \rightarrow \R \cup \{\infty\}$, $\dom(\psi) \coloneqq  \{x \in \R^q \mid \psi(x) < \infty\}$;  its subdifferential operator is
	$\partial \psi: \dom(\psi) \rightrightarrows \R^q:x\mapsto  \{ v \in \R^q \mid \psi(z) \geq \psi(x) + \langle v \mid z-x \rangle , \forall  z \in {\rm dom}(\psi) \}$; if $\psi$ is differentiable and convex, $\partial \psi=\nabla \psi$. For a set  $S\subseteq \R^q$, $\iota_{S}:\R^q \rightarrow \{ 0, \infty \}$ is the indicator function, i.e., $\iota_{S}(x) = 0$ if $x \in S$, $\infty$ otherwise; $\nc_{S}: S \rightrightarrows \R^q:x\mapsto \{ v \in \R^q \mid \sup_{z \in S} \, \langle v \mid z-x \rangle \leq 0  \}$ is the normal cone operator of $S$. If $S$ is closed and convex, then $\partial \iota_S=\nc_S$
	 and $(\Id+\nc_S)^{-1}=\proj_S$ is the Euclidean projection onto  $S$.
Given $\mc{F}:S\rightarrow \R^q$, the variational inequality VI$(\Fc,S)$ is the problem of finding $x^*\in S$ such that $\langle \Fc(x^*)\mid x-x^*\rangle \geq 0$, for all $x \in S$ (or, equivalently, $x^*$ such that $\0\in\Fc(x^*)+\nc_S(x^*)$). We denote  the solution set of VI$(\Fc,S)$ by SOL$(\Fc,S)$.

	\section{Mathematical setup}\label{sec:mathbackground}
	\begin{algorithm*} \caption{  Fully-distributed v-GNE seeking via \gls{PPPA} } \label{algo:1}
		\vspace{0.4em}
		\begin{itemize}[leftmargin=6.3em]
			\item[Initialization:]
			\begin{itemize}[leftmargin=1em]
				\item For all $i\in \mc{I}$, set $x_i^0\in \Omega_i$, $\bs{x}_{i,-i}^0\in \R^{n-n_i}$, $z_i^0=\0_m$, $\lambda_{i}^0\in \R^m_{\geq 0}$.
				\end{itemize}
			\vspace{0.1em}
			\item[For all $\k$:]
			\begin{itemize}[leftmargin=1em]
				\item
				Communication: The  agents exchange the variables  $\{ x^k_{i},\bs{x}_{i,-i}^k,\lambda_i^k \}$ with their neighbors.
				\item Local variables update: each agent $i\in \mc{I}$  computes
				\begin{align*}
				\xbs_{i,-i}^{k+1}   & = \textstyle \frac{1}{1+ \tau_i d_i}( {\xbs}_{i,-i}^k+\tau_i   \sum_{j \in \mc{N}_i} w_{i,j} \bs x_{j,-i}^k)
				\\
				x_i^{k+1}         & =
				\underset{y \in \Omega_i} {\argmin} \bigl(
				J_i(y,\xbs_{i,-i}^{k+1})+\textstyle \frac{1}{2\alpha\tau_i }  \bigl\|  y-x_i^{k} \bigr\|^2
				+\textstyle \frac{d_i}{2\alpha }
				\bigl\| \textstyle   y- \frac{1}{d_i}\sum_{j \in \mc{N}_i} w_{i,j}  \xbs_{j,i}^k \bigr\|^2+\frac{1}{\alpha}(A_i^\top \lambda_i^k)^\top y  \bigr)
				\\
				{z_i}^{k+1}       & = z_i^{k}+\textstyle  \sum_{j \in \mc{N}_i} \nu_{(i,j)} w_{i,j}  (\lambda_i^k-\lambda_{j}^k)
				\\
				{\lambda}_i^{k+1} & =\proj_{\R^{m}_{\geq  0}}
				\bigl(
				{\lambda}_i^k +\textstyle{\delta_i}
				\bigl( A_i(2x_i^{k+1}-x_i^k)-b_i-(2z_i^{k+1}-z_i^k)
				\bigr)
				\bigr).
				\end{align*}
			\end{itemize}
		\end{itemize}
		\vspace{-0.6em}
	\end{algorithm*}
	We consider a set of  agents, $ \mc I\coloneqq \{ 1,\ldots,N \}$, where each agent $i\in \mc{I}$ shall choose its decision variable (i.e., strategy) $x_i$ from its local decision set $\textstyle \Omega_i \subseteq \R^{n_i}$. Let $x \coloneqq  \col( (x_i)_{i \in \mc I})  \in \Omega $ denote the stacked vector of all the agents' decisions, $\textstyle \Omega \coloneqq  \Omega_1\times\dots\times\Omega_N\subseteq \R^n$ the overall action space and $\textstyle n\coloneqq \sum_{i=1}^N n_i$.
	The goal of each agent $i \in \mc I$ is to minimize its objective function $J_i(x_i,x_{-i})$, which depends on both the local variable $x_i$ and on the decision variables of the other agents $x_{-i}\coloneqq  \col( (x_j)_{j\in \mc I\backslash \{ i \} })$.
	Furthermore, the feasible decisions of each agent depends   on the action of the other agents via coupling constraints, which we assume affine: most of the literature focuses on this case \cite{PariseGentileLygeros_TCNS2020}, \cite{BelgioiosoNedicGrammatico2020},  which in fact accounts for the vast  majority of practical applications \cite[§3.2]{FacchineiKanzow2010}. Specifically, the overall feasible set is
	%
	$\mc{X} \coloneqq  \Omega \cap\left\{x \in \R^{n} \mid Ax\leq b \right\}$,
	where $A\coloneqq \left[A_{1}, \ldots, A_{N}\right]$ and $b\coloneqq \sum_{i=1}^{N} b_{i}$,  $A_{i} \in \R^{m \times n_{i}}$ and $b_{i} \in \R^{m}$ being local data. The game is then represented by the inter-dependent optimization problems:
	\begin{align} \label{eq:game}
	\forall i \in \mc{I}:
	\  \minimize{ y_i \in \R^{ {n_{i}}}}   \; J_i(y_i,x_{-i}) \quad
	\text{s.t.}  \  (y_i,x_{-i}) \in \mc X.
	\end{align}
	The technical problem we consider here is the computation of a \gls{GNE}, namely a set of decisions that simultaneously solve all the optimization problems in  \eqref{eq:game}.
	\begin{defn1}
		A collective strategy $x^{*}=\operatorname{col}\left((x_{i}^{*}\right)_{i \in \mathcal{I}})$ is a generalized Nash equilibrium if, for all $i \in \mc{I}$,
	$
		J_{i}\left(x_i^*, x_{-i}^{*}\right)\leq \inf \{J_{i}\left(y_{i}, x_{-i}^{*}\right) \mid (y_i,x_{-i}^*)\in\mc{X} \}.
		$ \hfill $\square$
	\end{defn1}
	Next, we postulate some common regularity and convexity
	assumptions for the constraint sets and cost functions, as in, e.g., \cite[Asm.~1]{Koshal_Nedic_Shanbag_2016}, \cite[Asm.~1]{Pavel2018}.
	\begin{standing}\label{Ass:Convexity}
		For each $i\in \mathcal{I}$, the set $\Omega_i$ is closed and convex; $\mc{X}$ is non-empty
		and satisfies Slater's constraint qualification;  $J_{i}$ is continuous and $J_{i}\left(\cdot, x_{-i}\right)$ is convex and continuously differentiable for every $x_{-i}$.
		{\hfill $\square$} \end{standing}
	As per standard practice \cite{PariseGentileLygeros_TCNS2020}, \cite{YiPavel2019}, among all the possible \glspl{GNE}, we focus on the subclass of \glspl{v-GNE} \cite[Def.~3.11]{FacchineiKanzow2010}, which are more economically justifiable, as well as computationally tractable  \cite{KulkarniShanbag2012}.
 The \glspl{v-GNE} are so called because they coincide with the solutions to the variational inequality VI$(F,\mc{X})$,
	where $F$ is the \emph{pseudo-gradient} mapping of the game:
	\begin{align}
	\label{eq:pseudo-gradient}
	F(x)\coloneqq \operatorname{col}\left( (\nabla _{\!x_i} J_i(x_i,x_{-i}))_{i\in\mathcal{I}}\right).
	\end{align}
	Under Standing Assumption~\ref{Ass:Convexity}, $x^*$ is a \gls{v-GNE} of the game in \eqref{eq:game} if and only if there exists a dual variable $\lambda^*\in \R^m $ such that the following \gls{KKT} conditions are satisfied \cite[Th.~4.8]{FacchineiKanzow2010}:
	\begin{align} \label{eq:KKT}
	\begin{aligned}{\0_{n}} & {\in F\left(x^{*}\right)+A^\top \lambda^{*}+\mathrm{N}_{\Omega}\left(x^{*}\right)}
	\\
	{\0_{m}}                  & {\in-\left(A x^{*}-b\right)+\mathrm{N}_{\R_{\geq 0}^{m}}\left(\lambda^{*}\right)}.\end{aligned}
	\end{align}
	\begin{standing}\label{Ass:StrMon}
		The pseudo-gradient mapping $F$ in \eqref{eq:pseudo-gradient}  is $\mu$-strongly monotone and $\theta_0$-Lipschitz continuous, for some $\mu$, $\theta_{0}>0$.
		\hfill $\square$
	\end{standing}
The strong monotonicity of $F$  is sufficient to ensure existence and uniqueness of a \gls{v-GNE} \cite[Th. 2.3.3]{FacchineiPang2007}; it was always assumed for \gls{GNE} seeking under partial-decision information with fixed step sizes \cite[Asm.~2]{TatarenkoShiNedicTAC20}, \cite[Asm.~3]{Pavel2018} (while it is sometimes replaced by strict monotonicity or cocoercivity, under vanishing steps and compactness of $\mc{X}$
\cite[Asm.~2]{Koshal_Nedic_Shanbag_2016}, \cite[Asm.~3]{PangHu_TAC2020} \cite[Asm.~5]{BelgioiosoNedicGrammatico2020}).

	\section{Fully-distributed  equilibrium seeking }\label{sec:distributedGNE}

	In this section, we present our baseline algorithm to seek a \gls{v-GNE} of the game in  \eqref{eq:game} in a fully-distributed way.
	Specifically,
	each agent $i$ only knows its own cost function $J_i$ and feasible set $\Omega_i$, and the portion of the coupling constraints  $(A_i,b_i)$.
	Moreover, agent $i$ does not have full knowledge of $x_{-i}$, and only relies on the information exchanged locally with some neighbors over an undirected  communication network $\mathcal G(\mc{I},\mc{E})$.  The unordered pair $(i,j) $ belongs to the set of edges $\mc{E}$ if and only if agent $i$ and $j$ can mutually exchange information.
	We denote: $W=[w_{i,j}]_{i,j\in \mc{I}}\in \R^{N\times N}$ the  symmetric weight matrix of $\mc{G}$, with $w_{i,j}>0$ if $(i,j)\in \mc{E}$, $w_{i,j}=0$ otherwise, and the convention $w_{ii}=0$ for all $i\in\mc{I}$;
	$L\coloneqq D-W$ the Laplacian matrix of $\mc{G}$, with degree matrix $D\coloneqq \diag((d_i)_{i\in\mc{I}})$, and $d_{i}\coloneqq \textstyle \sum_{j=1} ^N w_{i,j}$ for all $i\in \mc{I}$; $\mc{N}_i=\{j\mid (i,j)\in \mc{E}\}$ the set of neighbors of agent $i$. Moreover, we label the edges $(e_\ell)_{\ell\in\{1,\dots,E\}}$, where $E$ is the cardinality of $\mc{E}$, and we assign to each edge $e_\ell$ an arbitrary orientation.
	We denote  the weighted incidence matrix as $V\in \R^{E\times N}$, where $[V]_{\ell,i}=\sqrt{(w_{i,j})}$ if $e_\ell=(i,j)$ and $i$ is the output vertex of $e_\ell$, $[V]_{\ell,i}=-\sqrt{(w_{i,j})}$ if $e_\ell=(i,j)$ and $i$ is the input vertex of $e_\ell$, $[V]_{\ell,i}=0$ otherwise.
	It holds that $L=V^\top V$; moreover, $\Null(V)=\Null(L)=  \{\kappa\1_N,\kappa\in\R\}$ under the following connectedness assumption \cite[Ch.~8]{Godsil:algebraic_graph_theory}.


	\begin{standing}
		\label{Ass:Graph}
		The communication graph $\mathcal G (\mc{I},\mc{E}) $ is undirected and connected.
		\hfill $\square$
	\end{standing}
	%
	%
	%
	%
	In the partial-decision information, to cope with the lack of knowledge, each agent keeps an estimate of all other agents' actions \cite{YeHu2017}, \cite{TatarenkoShiNedicTAC20}, \cite{Pavel2018}. We denote $\xbs_{i}\coloneqq \operatorname{col}((\xbs_{i,j})_{j\in \mc{I}})\in \R^{n}$,  where $\xbs_{i,i}\coloneqq x_i$ and $\xbs_{i,j}$ is agent $i$'s estimate of agent $j$'s action, for all $j\neq i$; let also $\xbs_{j,-i}\coloneqq \col((\xbs_{j,\ell})_{\ell\in\mc{I}\backslash\{i\}})$. Moreover, we let each agent  keep an estimate $\lambda_i\in \R_{\geq 0}^m$ of the dual variable, and an auxiliary variable $z_i\in\R^m$.

	Our proposed dynamics are summarized in Algorithm~\ref{algo:1}, where the global parameter $\alpha>0$ and the positive step sizes $\tau_i$, $\delta_i$, $\nu_{(i,j)}=\nu_{(j,i)}$, for all $i \in \mc{I}$ and $(i,j)\in\mc{E}$,  have to be chosen appropriately (see §\ref{sec:derivationconvergence}).
	Each agent $i$ updates its action $x_i$ similarly to a proximal best-response, but with two extra terms that are meant to penalize and correct the disagreement among the estimates and the coupling constraints violation. Most importantly, the agents
	evaluate their cost functions in their local estimates, not on the actual collective strategy. In steady state, the agents should agree on their estimates, i.e., $\xbs_i=\xbs_j$,  $\lambda_i=\lambda_j$, for all $i,j\in\mc{I}$.	This motivates the presence of consensus terms for both primal and dual variables.
	From a control-theoretic perspective,
 the updates of each $z_i$   can be seen  as integrator dynamics driven by the disagreement of the variables $\lambda_j$’s. This integral action is meant to permit the distributed asymptotic satisfaction of the coupling constraints, despite the computation of each $\lambda_i$ only involves the local block $(A_i, b_i)$ -- differently from typical centralized dual ascent iterations.  We postpone a formal derivation of Algorithm~\ref{algo:1} to §\ref{sec:derivationconvergence}.
	\begin{rem1}\label{rem:1_singlevaluedargmin}
		The functions $J_i(\cdot,\xbs_{i,-i})$ are strongly convex, for all $\xbs_{i,-i}$, $i\in \mc{I}$, as a consequence of Standing Assumption~\ref{Ass:StrMon}. Hence, the $\argmin$ operator in Algorithm~\ref{algo:1} is single-valued, and the algorithm is well defined.
		\hfill $\square$
		\end{rem1}
	     \begin{rem1}
	     	In Algorithm~\ref{algo:1}, each agent has to  \emph{locally} solve an optimization problem, at every iteration. Not only these subproblems are
	     	fully-decentralized (i.e., they do not require extra communication), but they are also of low dimension ($n_i$). This is a major departure from the procedure proposed in the \glspl{PPA} \cite[Alg.~2]{Scutari_ComplexGames_TIT2014}, \cite[Alg.~2]{YiPavel_PPP_2019},
		 where the agents have to collaboratively solve a subgame (of dimension $n$) before each update.  \hfill $\square$
	\end{rem1}

	\section{Convergence analysis }\label{sec:derivationconvergence}
	\subsection{Definitions and preliminary results}
	 We denote $\xbs\coloneqq \col((\xbs_i)_{i\in\mc{I}})\in\R^{Nn}$.
	  Besides, let us define, as in \cite[Eq.~13, 14]{Pavel2018},
	   for all $i \in \mc{I}$,
	\begin{subequations}
		\begin{align}
		\mathcal{R}_{i}\coloneqq & \begin{bmatrix}{{\0}_{n_{i} \times n_{<i}}} & {I_{n_{i}}} & {\0_{n_{i} \times n_{>i}}}\end{bmatrix} \in \R^{n_i \times n },
		\\
		\mathcal{S}_{i}\coloneqq  &\begin{bmatrix}{I_{n_{<i}}} & {\0_{n_{<i} \times n_{i}}} & {\0_{n_{<i} \times n_{>i}}} \\ {\0_{n_{>i} \times n_{<i}}} & {\0_{n_{>i} \times n_{i}}} & {I_{n_{>i}}}\end{bmatrix} \in \R^{n_{-i}\times n }\hspace{-2pt}
		\end{align}
	\end{subequations}
	where $n_{<i}\coloneqq \sum_{j<i,j \in \mathcal{I}} n_{j}$, $n_{>i}\coloneqq \sum_{j>i, j \in \mathcal{I}} n_{j}$ and $n_{-i}\coloneqq n-n_i$. In simple terms, $\mathcal R _i$ selects the $i$-th $n_i$-dimensional component from an $n$-dimensional vector, while $\mathcal S_i$ removes it. Thus, $\mathcal{R}_{i} \xbs_{i}=\xbs_{i,i}=x_i$ and $\mathcal{S}_{i} \xbs_{i}=\xbs_{i,-i}$. Let $\mathcal{R}\coloneqq \operatorname{diag}\left((\mathcal{R}_{i})_{i \in \mathcal{I}}\right)$,  $\mathcal{S}\coloneqq \operatorname{diag}\left((\mathcal{S}_{i})_{i \in \mathcal{I}}\right)$. It follows that $x=\mathcal{R} \xbs$ and $\operatorname{col}((\xbs_{i,-i})_{i \in \mathcal{I}})=\mathcal{S}{\xbs} \in \R^{(N-1) n}$. Moreover,
	$\xbs=\mathcal{R}^\top x+\mathcal{S}^\top \mc{S} \xbs.$
	We define the  \textit{extended pseudo-gradient} mapping  $\bs{F}:\R^{Nn}\rightarrow \R^n$ as
	\begin{align}
	\label{eq:extended_pseudo-gradient}
	\bs{F}(\xbs)\coloneqq \operatorname{col}\left((\nabla_{\!x_{i}} J_{i}\left(x_{i}, \xbs_{i,-i}\right))_{i \in \mathcal{I}}\right),
	\end{align}
	and the operators
	\begin{align}
	\label{eq:Fa}
	\Fa(\xbs) & \coloneqq \alpha\mc{R}^\top\bs{F}(\xbs)+(\bs{D}_{\!n}-\bs{W}_{\!\!n})\xbs,
	\\
	\label{eq:opA}
	\mathcal{A}(\obs)                & \coloneqq  \hspace{-0.2em}
	\underbrace{\left[
	\begin{matrix}
	\Fa(\bs{x})\\
	\0_{Em}\\
	\bs b
	\end{matrix}\right]
	\hspace{-0.2em}+\hspace{-0.2em}\left[
	\begin{matrix}
	{\mathcal{R}^\top\bs{A}^\top \bs{\lambda}}  \\ {-\bs{V}_{\!\!m} \bs{\lambda}} \\ {{\bs{V}_{\!\!m}^{\top} \bs{v} }  -\bs{A} \mathcal{R}}\xbs
	\end{matrix}\right]}_{\coloneqq \Ac_1(\omega)}
	\hspace{-0.2em}+\hspace{-0.2em}
	\left[
	\begin{matrix}
	\nc_{\bs{\Omega}}(\xbs) \\
	\0_{Em}\\
	\nc_{\R^{Nm}_{\geq 0}}(\bs{\lambda})
	\end{matrix}\right]\hspace{-0.5 em}
	\end{align}
	\sloppy
	where $\alpha>0$ is a design constant,
	$\obs\coloneqq \col(\xbs,\bs{v},\bs{\lambda})$,    $\bs{v}\coloneqq \col((v_{\ell})_{\ell\in \{1,\dots,E\}})\in\R^{Em}$, $\bs{\lambda}\coloneqq  \col((\lambda_i)_{i\in\mc{I}})\in\R^{Nm}$, $\bs{A}\coloneqq \diag((A_i)_{i\in\mc{I}})$,
	$\bs{W}_{\!\!n}\coloneqq W\otimes I_n$,
	$\bs{D}_{\!n}\coloneqq D\otimes I_n$,
	$\bs{V}_{\!\!m} \coloneqq V\otimes I_m$,
	and
	$\bs{\Omega}\coloneqq \{\xbs\in\R^{nN}\mid \mc{R}\xbs\in \Omega\}$.

	The following lemma relates the unique \gls{v-GNE} of the game in \eqref{eq:game} to the zeros of the operator $\mc{A}$. The proof is analogous to \cite[Th.~1]{Pavel2018} or Lemma~\ref{lem:zerosagg} in  Appendix~\ref{app:th:main1agg},  and hence it is omitted.
	\begin{lem1}\label{lem:zeros}
			Let $\Ac$ be as in \eqref{eq:opA}. It holds that $\zer(\Ac)\neq \varnothing$. Moreover, let $\xbs^*\in \R^{Nn}$, $\bs{\lambda}^*\in\R^{Nm}$; then, the following statements are equivalent:
			\begin{lemlist}[topsep=-2em]
				\item\label{lem:zerosi} There exists $\vbs^*$ such that $\col(\xbs^*,\bs{v}^*,\bs{\lambda}^*)\in\zer(\mc{A})$.
				\item\label{lem:zerosii} $\xbs^*=\1_N\otimes x^*$ and $\lbs^*=\1_N\otimes \l^*$, where the pair $(x^*,\l^*)$ satisfies the \gls{KKT} conditions in \eqref{eq:KKT}, hence $x^*$ is the \gls{v-GNE} of the game in \eqref{eq:game}.  	\hfill $\square$
			\end{lemlist}
		\end{lem1}
	 Effectively, Lemma~\ref{lem:zeros} provides an extension of the \gls{KKT} conditions in \eqref{eq:KKT} and allows us to recast the \gls{GNE} problem as that of computing a zero of the operator $\Ac$, for which a number of iterative algorithms are available \cite[§26-28]{Bauschke2017}.
In fact, in §\ref{subsubsec:derivation}, we show  that  Algorithm~\ref{algo:1} can be recast as a \gls{PPA} \cite[Th.~23.41]{Bauschke2017}.

	Nonetheless, technical difficulties arise in the analysis because of the  partial-decision information setup.
	 Specifically, in  \eqref{eq:extended_pseudo-gradient},  each partial gradient $\nabla_{\!x_{i}} J_{i}\left(x_{i}, \xbs_{i,-i}\right)$ is evaluated on the local estimate $\xbs_{i,-i}$, and not on the actual value $x_{-i}$. Only  when the estimates $\bs{x}$ are at consensus, i.e., $\bs{x}=\1_N\otimes x$ (namely, the estimate of each agents coincide with the actual value of $x$), we have that $\bs{F}(\bs{x})=F(x)$. As a result, the operator $\mc{R}^\top \bs{F}$ (and consequently the operator $\Ac$) is not monotone in general%
	 \footnote{It can be shown that $\mc{R}^\top \bs{F}$ is monotone only if the mappings $\nabla_{ \!x_{i}} J_{i}(x)$'s do not depend on  $x_{-i}$ (in which case, there is no need for a partial-decision information assumption).}\!,
	  not even under the strong monotonicity of the game mapping $F$ in Standing Assumption~\ref{Ass:StrMon}.
	 Instead,  analogously to the approaches in \cite{SalehisadaghianiWeiPavel2019}, \cite{Pavel2018}, \cite{GadjovPavel2018}, our analysis is based on a \emph{restricted} monotonicity property.
 \begin{defn1}\label{def:restricted}
		An operator $\mathcal{F} : \R^q \rightrightarrows \R^q$ is
		restricted ($\mu$-strongly) monotone in $\mc{H}_P$ if $\zer(\Fc)\neq \varnothing$ and  $\langle \omega-\omega^*\mid u \rangle_P \geq 0$ ($\geq \mu\|\omega-\omega^*\|^2_P$) for all $(\omega,u)\in\gra(\Fc)$, $\omega^*\in\zer(\Fc)$ (we omit the characterization ``in $\mc{H}_P$'' whenever  $P=I$). \hfill $\square$
	\end{defn1}
   Definition~\ref{def:restricted} differs from that in \cite[Lem.~3]{Pavel2018}, as we  only consider properties with respect to the zero set and we need to include set-valued operators.  The definition comprises the  nonemptiness of the zero set and it does not exclude an operator that is multi-valued on its zeros. The next  lemmas show that restricted monotonicity of $\Ac$ can be guaranteed for any game satisfying
	Standing Assumptions~\ref{Ass:Convexity}--\ref{Ass:Graph}, without additional hypotheses.

	\begin{lem1}[{\cite[Lemma 3]{Bianchi_ECC20_ctGNE}}]\label{lem:LipschitzExtPseudo}
		The  mapping $\bs{F}$ in \eqref{eq:extended_pseudo-gradient} is $\theta$-Lipschitz continuous, for some $\theta\in [\mu, \theta_0]$.
		{\hfill $\square$} \end{lem1}
	\begin{lem1}\label{lem:strongmon_constant}
		Let  $\ \alphamax \coloneqq   \frac{4\mu\uplambda_2(L)}{(\theta_0+\theta)^{2}+4\mu\theta}$,
		\begin{equation}
		\textstyle
		\label{eq:alphamax}
		{M}\coloneqq \alpha\begin{bmatrix}{\frac{\mu}{N}} & \ {-\frac{\theta_0+\theta}{2\sqrt{N}}} \\ {-\frac{\theta_0+\theta}{2\sqrt{N}}} & \ { \textstyle \frac{\uplambda_{2}(L)}{\alpha}-\theta}\end{bmatrix} \!\!, \
		~~~
		\begin{aligned}
		\mu_{\Fa}& \coloneqq \uplambda_{\textnormal{min}}(M).
		\end{aligned}
		\end{equation}
		If $\alpha\in (0,\alphamax]$ , then $\mu_{\Fa}\geq 0$ and $\mc{A}$ in \eqref{eq:opA} is restricted monotone. \hfill $\square$
	\end{lem1}
	\begin{proof1}
		The operator $\Ac$ in \eqref{eq:opA} is the sum of three operators. The third is monotone by properties of normal cones \cite[Th.~20.25]{Bauschke2017}; the second is a linear skew-symmetric operator, hence monotone \cite[Ex.~20.35]{Bauschke2017}. Let   $\obs^*=\col(\xbs^{*},\bs{v}^{*},\bs{\lambda}^{*})\in\zer(\mc{A})$, where $\zer(\mc{A})\neq \varnothing$ by Lemma~\ref{lem:zeros}. By Lemma~\ref{lem:zeros}, $\xbs^*=\1_N \otimes x^*$, with $x^*$ the \gls{v-GNE} of the game in \eqref{eq:game}; hence by \cite[Lemma~3]{Pavel2018}, for any $\alpha\in (0,\alphamax]$, it holds that $M\succcurlyeq0$ and that, for all $\xbs\in\R^{Nn}$
		\begin{equation}\label{eq:strongFa}
		\langle \xbs-\xbs^{*}
		\mid \Fa (\xbs)-\Fa (\xbs^*)\rangle  \geq \mu_{\Fa}\|\xbs-\xbs^*\|^2.
		\end{equation}
	Therefore, for all $(\obs,\bs{u})\in\gra(\mc{A})$, with $\obs=\col(\xbs,\bs{v},\bs{\lambda})$, $\langle \obs-\obs^* \mid \bs{u}-\0 \rangle \geq  \mu_{\Fa} \|\xbs-\xbs^*\|^2 \geq 0$. \hfill $\blacksquare$
	\end{proof1}
	\subsection{\gls{PPA} for restricted monotone operators}\label{subsec:PPA}
	In the remainder of this section, we show that Algorithm~\ref{algo:1} is an instance of the \gls{PPA}, applied to seek a zero of the (suitably preconditioned) operator $\mc{A}$ in \eqref{eq:opA}. Then, we show its convergence based on the restricted monotonicity result  in Lemma~\ref{lem:strongmon_constant}.

	Informally speaking, in proximal-point methods,	 a problem is decomposed into a sequence of regularized subproblems, which are  possibly better conditioned and  easier to solve.
	Let $\mc{B}:\R^q\rightrightarrows \R^q$ be  maximally monotone   \cite[Def.~20.20]{Bauschke2017} in a space $\mc{H}_P$, and $\J_\mc{B}=(\Id+\mc{B})^{-1}$ its resolvent. Then, $\dom (\J_\mc{B})=\R^q$ and $\J_\mc{B}$ is single-valued;  moreover, if
	 $\zer({\mc{B}})\neq \varnothing$, then the sequence $(\omega^k)_\k$ generated by the \gls{PPA},
	\begin{equation}\label{eq:PPPmaximally}
	(\forall k\in\N) \quad \omega^{k+1}= \J_\mc{B}(\omega^{k}), \quad \omega^0\in\R^q,
	\end{equation}
	converges to a point in $\zer(\mc{B})=\fix(\J_\mc{B})$ \cite[Th.~23.41]{Bauschke2017}. Note that   performing the update in \eqref{eq:PPPmaximally} is equivalent to solving for $\omega^{k+1}$ the (regularized) inclusion
	\begin{equation}\label{eq:inclusionformulation}
	\0\in \mc{B}(\omega^{k+1}) +\omega^{k+1}-\omega^k.
	\end{equation}
	  Unfortunately, many operator-theoretic properties are not guaranteed if $\mc{B}$ is only restricted monotone.
       In fact, $\J_\mc{B}$ might not be defined everywhere
	or single-valued.
\begin{exmp1}
		Let $\mc{B}:\R\rightarrow \R$, with  $\mc{B}(\omega)=9-2\omega$ if  $\omega\in [3,4)$,
		 $\mc{B}(\omega)=\omega$ otherwise. Then, $\zer(\mc{B})=\{0\}$  and  $\mc{B}$ is restricted strongly monotone. However, $\J_{\mc{B}}(\omega)=\{\frac{\omega}{2},9-\omega\}$ if $\omega\in[5,6)$ and  $\J_{\mc{B}}(\omega)=\varnothing$ if $\omega\in(6,8)$. \hfill $\square$
	\end{exmp1}
Nonetheless, some  important properties carry on to the restricted monotone case, as we prove next.
\begin{lem1}\label{lem:FQNE}
Let $\mc{B}:\R^q\rightrightarrows \R^q $ be restricted monotone in $\mc{H}_P$.
Then, $\J_{\mc{B}}$ is firmly quasinonexpansive in $\mc{H}_P$: for any $(\omega,u)\in\gra(\J_{\mc{B}})$, $\omega^*\in \zer(\mc{B})=\fix(\J_{\mc{B}})$,
 it holds that
	\begin{equation}\label{eq:FQNE}
	   \langle \omega-u\mid \omega-\omega^*\rangle_P-\|u-\omega\|_P^2 =\langle \omega-u \mid u-\omega^*\rangle_P \geq 0.
	\end{equation}
 Moreover, $\J_{\mc{B}}(\omega^*)=\{ \omega^*\} $. \hfill$\square$
\end{lem1}
\begin{proof1}
By definition of resolvent, $ \omega^*\in\J_{\mc{B}} (\omega^*) \iff \omega^*+\mc{B}\omega^* \ni \omega^*\iff \0 \in  \mc{B}(\omega^*)$;
also, for any $(\omega,u)\in\gra(\J_{\mc{B}})$, $\omega-u\in \mc{B}(u)$. Hence, the inequality in \eqref{eq:FQNE} is the restricted monotonicity of $\mc{B}$; the elementary equality follows by expanding the terms. Finally, by taking $\omega=\omega^*$ in \eqref{eq:FQNE}, we infer that $\J_{\mc{B}}$ is single-valued on $\fix(\J_{\mc{B}})$. \hfill$\blacksquare$
\end{proof1}
Next, by leveraging Lemma~\ref{lem:FQNE}, we extend classical results for the \gls{PPA} \cite[Th.~5.6]{Combettes2001} to the case of a restricted monotone  operator (possibly with multi-valued resolvent).
 \begin{thm}\label{th:FNEconvergence}
 	Let $\mc{B}:\R^q\rightrightarrows \R^q$ be restricted monotone in $\mc{H}_P$, and  $C\coloneqq \zer(\mc{B})\neq \varnothing$.
     Let $(\gamma^k)_{k\in\N}$ be a sequence in $[0,2]$, and $(e^k)_{k\in\N}$ a sequence in $\R^q$ such that $(\gamma^k\|e^k\|_P)_{k\in\N} \in \ell^1$. Let $\omega^0\in\R^q$ and  let  $(\omega^k)_{k\in\N}$ be any sequence such that:
	\begin{equation}\label{eq:KM}
	(\forall k\in\N) \quad \omega^{k+1} = \omega^k+\gamma^k(u^k-\omega^k+e^k),  u^k\in \J_{\mc{B}}(\omega^k).
	\end{equation}
	Then, the following statements hold:
	\begin{thmlist}
		\item \label{th:FNEconvergencei}$(\forall \omega^*\in C)(\forall k\in\N) \;\;\; \|\omega^{k+1}-\omega^*\|_P\leq \|\omega^k-\omega^*\|_P+\gamma^k\|e^k\|_P $.
		\item \label{th:FNEconvergenceii}$ \left(\gamma^k(2-\gamma^k)\|u^k-\omega^k\|_P^2\right)_{k\in\N}\in \ell^1$.
		\item \label{th:FNEconvergenceiii}  Assume that every cluster point of $(\omega^k)_{k\in\N}$ belongs to $C$. Then, $(\omega^k)_{k\in\N}$  converges to a point in $C$.
		\item\label{th:FNEconvergenceiv}Assume that   $\mc{B}$ is $\mu_{\mc{B}}$-strongly restricted
		 monotone in $\mc{H}_P$. Then,   $C=\{ \omega^* \}$ and,
		 for all $\k$, $\|\omega^{k+1}-\omega^*\|_P\leq \rho^k \|\omega^{k}-\omega^*\|_P+ \gamma^k\|e^k\|_P$,
		  where  $\rho^k=\max(1-\frac{\gamma^k \mu_{\mc{B}}}{1+\mu_{\mc{B}}},\gamma^k-1)$. \hfill$\square$
	\end{thmlist}
\end{thm}
\begin{proof1}
	See Appendix~\ref{app:th:FNEconvergence}. \hfill $\blacksquare$
	\end{proof1}
	\begin{rem1}
	The condition $\dom(\J_{\mc{B}})=\R^q$  is sufficient (but not necessary) for
    the existence of a sequence $(\omega^k)_{\k}$ that satisfies \eqref{eq:KM}, which can be constructed  choosing arbitrarily $u_k\in\J_{\mc{B}}(\o^k)$, for all $\k$. \hfill $\square$
    \end{rem1}
\begin{exmp1}\label{ex:pseudomonotone}
	Consider the \textnormal{VI}$(\Psi,S)$, where $S\subset\R^q$ is compact and convex, and $\Psi:\R^q\rightarrow \R^q$ is continuous and pseudomonotone in the sense of Karamardian (i.e., for all $\omega,\omega^\prime \in\R^q$, the  implication $\langle \Psi(\omega),\omega^\prime-\omega \rangle \geq 0  \Rightarrow \langle \Psi(\omega^\prime),\omega^\prime-\omega \rangle \geq 0
	$ holds). It holds that SOL$(\Psi,S)=\zer(\mc{B}) \neq \varnothing$ , where $\mc{B}=\Psi+\nc_{S}$  \cite[Prop.~2.2.3]{FacchineiPang2007}. Moreover $\mc{B}$ is restricted monotone. To show this, consider any $\omega^*\in\zer(\mc{B})$ and $(\omega,u)\in\gra(\mc{B})$, so  $u=\Psi(\omega)+u^\prime$, for some $u^\prime$ such that $(\omega,u^\prime)\in\gra(\nc_{S})$. Then, $\langle u \mid \omega-\omega^*\rangle=\langle \Psi(\omega) \mid \omega-\omega^*\rangle+\langle u^\prime-\0 \mid \omega-\omega^*\rangle\geq 0$,
	where we used that $\langle \Psi(\omega) \mid \omega-\omega^*\rangle\geq 0$, by pseudomonotonicity and because  $\langle \Psi(\omega^*) \mid \omega-\omega^*\rangle\geq 0$ by definition of \textnormal{VI},
	 and $\langle u^\prime-\0 \mid \omega-\omega^*\rangle\geq 0$
	 because $(\omega^*,\0)\in\gra(\nc_S)$ and  monotonicity of the normal cone.

	We note that $\dom(\J_{\mc{B}})=\R^q$ by \cite[Prop.~2.2.3]{FacchineiPang2007}. Let us consider any sequence $(\omega^k)_\k$ such that, for all $\k$, $\omega^{k+1}=u^k, u^k\in\J_{\mc{B}}(\omega^k)$, (or equivalently \eqref{eq:inclusionformulation} or  $\omega^{k+1}\in \textnormal {SOL}(\Psi+\Id-\omega^k,S)$).
	By Theorem~\ref{th:FNEconvergence} (with $\gamma^k=1, e^k=\0$), $(\omega^k)_\k$ is bounded, hence it admits at least one cluster point, say $\bar{\omega}$; by  Theorem~\ref{th:FNEconvergenceii} $\|u^k-\omega^k\|\rightarrow 0$. However, by definition of  VI, for any $\omega \in S$, $\langle \Psi(u^k)+u^k-\omega^k \mid \omega-u^k\rangle \geq 0 $.
	By passing to the limit (on a subsequence) and by continuity, we obtain
	$\langle \Psi(\bar{\omega})\mid \omega-\bar{\omega }\rangle \geq 0$,
	which shows that $\bar{\omega}\in \textnormal{SOL}(\Psi,S)$. Therefore $(\omega^k)_{\k}$ converges to a solution to  \textnormal{VI}$(\Psi,S)$ by  Theorem~\ref{th:FNEconvergenceiii}. This extends the results in \cite[§4.2]{ElFarouq_PseudomonotonePPP_2001}, where  hypomonotonicity of $\Psi$ is assumed and where a small-enough step size is chosen to ensure that $\J_{\mc{B}}$ is single-valued (besides, pseudomonotonicity of $\Psi$ is  sufficient, but not necessary, for the restricted monotonicity of $\mc{B}$, and Theorem~\ref{th:FNEconvergence} would also allow to take into account iterations with errors, cf. \cite[§4.2]{ElFarouq_PseudomonotonePPP_2001}). \hfill $\square$
\end{exmp1}

\subsection{Derivation and convergence}\label{subsubsec:derivation}
   Next, we show how that Algorithm~\ref{algo:1} is  obtained by applying the iteration in  \eqref{eq:KM} to the operator
	$\Phi^{-1} \mc{A}$, where
	\begin{equation}\label{eq:Phi}
	\Phi\coloneqq  \left[
	\begin{matrix}
	{\bar{\tau}^{-1}+\bs{W}_{\!\!n}} & {\0} & -{\mathcal{R}^\top\bs{A}^\top} \\ {\0} & {\bar{\nu}^{-1}} & {\bs{V}_{\!\!m} } \\ {-\bs{A} \mathcal{R}} & {\bs{V}_{\!\!m}^\top } & {\bar{\delta}^{-1}}
	\end{matrix}\right]
	\end{equation}
	\sloppy
	is called \emph{preconditioning} matrix. The step sizes $\bar{\tau}\coloneqq \diag((\tau_i I_n)_{i\in\mc{I}})$, $\bar{\nu}\coloneqq \diag((\nu_{(i,j)} I_m)_{(i,j)\in\mc{E}}) $, $\bar{\delta}\coloneqq \diag((\delta_i I_ m)_{i\in\mc{I}})$, have to be chosen such that  $\Phi\succ 0$. In this case, it also holds that
	$\zer(\Phi^{-1}\mc{A})=\zer(\mc{A})$. Sufficient conditions that ensure $\Phi\succ 0$ are given in the next lemma, which   follows by the   Gershgorin's circle theorem.
	\begin{lem1}\label{lem:stepsizes}
		The matrix $\Phi$ in \eqref{eq:Phi} is positive definite if ${\nu_{(i,j)}}^{-1} >  2\sqrt{(w_{i,j})}$ for all $(i,j)\in\mc{E}$  and $\tau_i ^{-1}> d_{i} +\|A_i^\top\|_{\infty}$, $\delta_i  ^{-1}> \|A_i \|_\infty+\textstyle \sum^N_{j=1} \sqrt{(w_{i,j})} $ for all $i\in\mc{I}$. \hfill $\square$
	\end{lem1}

	In the following, we always assume that the step sizes in Algorithm~\ref{algo:1} are chosen
	such that $\Phi\succ 0$. Then, we are able to formulate the following result.

	\begin{lem1}\label{lem:derivation}
		Algorithm~\ref{algo:1} is equivalent to the iteration
		\begin{equation}\label{eq:algcompact}
		(\forall k\in \mathbb{N}) \quad  \obs^{k+1} \in \J_{\Phi^{-1}\Ac}(\obs^k),
		\end{equation}
		with $\Ac$ as in \eqref{eq:opA}, $\Phi$ as in \eqref{eq:Phi}:  for any initial condition
		$\obs^0=\col(\xbs^0,\bs{v}^0=\0_{Em},\bs\lambda^0)$, the sequence $(\xbs^k,\bs{V}_{\!\!m}^\top\bs{v}^k,\bs\lambda^k)_{k\in\N}$ generated by \eqref{eq:algcompact} coincides with the sequence $(\xbs^k, \bs{z}^k,\bs\lambda^k)_{k\in\N}$ generated by Algorithm~\ref{algo:1} with initial conditions $(\xbs^0,\bs{z}^0=\0_{Nm},\bs\lambda^0)$.
		\hfill $\square$
	\end{lem1}
	\begin{proof1}
		By definition of inverse operator, we have that
		\allowdisplaybreaks
		\begin{align}
		\nonumber
		&
		\obs^{k+1} \in (\Id+\Phi^{-1}\Ac)^{-1}(\obs^{k})
		\\
		\nonumber
		\iff &
		\0  \in \Phi^{-1} \Ac(\obs^{k+1})-\obs^{k}+\obs^{k+1}
		\\
		\iff &
		\0  \in \Phi(\obs^{k+1}-\obs^{k})+\Ac(\obs^{k+1} ) \label{eq:explicitinclusion}
		\\
		\label{eq:useful1}
		\iff&
	\! \! \! \left\{
		\begin{aligned}
		\0 & \in
		 \bar{\tau}^{-1}(\xbs^{k+1}-\xbs^k)+\cancel{\bs{W}_{\! \! n} \xbs^{k+1}} -\bs{W}_{\!\!n}\xbs^k+\bs{D}_{\!n}\xbs^{k+1} \\
		  & \ha -\cancel{\mc{R}^{ \top}  \bs{A}^{\top }\bs{\lambda}^{k+1}}
		 +\mc{R}^{ \top} \bs{A}^{\top}\bs{\lambda}^{k} +\alpha\mc{R}^{ \top}\bs{F}(\xbs^{k+1})
		\\ & \ha   - \cancel{\bs{W}_{\!\!n}\xbs^{k+1}}+\cancel{\mc{R}^{\top} \bs{A}^{\top} \bs{\lambda}^{k+1} }+\nc_{\bs{\Omega}}(\xbs^{k+1})
		\\
		\0 & \in \bar{\nu}^{-1}(\bs{v}^{k+1}-\bs{v}^{k})+\cancel{\bs{V}_{\!\!m}\bs{\lambda}^{k+1}}-\bs{V}_{\!\!m}\bs{\lambda}^{k}-\cancel{\bs{V}_{\!\!m}\bs{\lambda}^{k+1}}
		\\
		\0 &
		\in \bar{\delta}^{-1}(\bs{\lambda}^{k+1}-\bs{\lambda}^k)+\nc_{\R^{mN}_{\geq  0}}(\bs{\lambda}^{k+1})+\bs{b} \\ & \ha -\bs{A}\mc{R}(2\xbs^{k+1}-\xbs^k)+\bs{V}_{\!\!m}^\top (2\bs{v}^{k+1}-\bs{v}^k)
		\end{aligned}
		\right.
		\raisetag{0.5\baselineskip}
		\end{align}
		In turn, the first inclusion in \eqref{eq:useful1} can be split in two  by left-multiplying both sides with $\mc{R}$ and $\mc{S}$. By $\mc{S}\mathrm{N}_{\bs{\Omega}}=\0_{(N-1)n}$, $\mc R\mc R^\top=I_n$ and  $\mc{S}\mc{R}^\top=\0_{(N-1)n\times n}$,  we get
		\begin{align}
		\nonumber
		&
		\left\{ \begin{aligned}
		\0                         & \in   \mc{S}((I+\bar{\tau}\bs{D}_{\!n})\xbs^{k+1}-\xbs^{k}-\bar{\tau}\bs{W}_{\!\!n}\xbs^{k})
		\\
		\0                              & \in 	\begin{aligned}[t] & \mc{R}((I+\bar{\tau}\bs{D}_{\!n})\xbs^{k+1}-\xbs^k-\bar{\tau}\bs{W}_{\!\!n}\xbs^k) \\ & +\mathrm{N}_{{\Omega}}(x^{k+1})
		+ \alpha\bar{\tau}\bs{F}((x^{k+1},\mc{S}\xbs^{k+1}))+\bar{\tau}\bs{A}^\top \bs{\lambda}^k
		\end{aligned}
		\end{aligned}
		\right.
		\\
		\nonumber
	\underset{\forall i\in \mc{I}}{\iff} \! &
		\left\{ \begin{aligned}
		\xbs_{i,-i}^{k+1}                      & = \textstyle \frac{1}{1+\tau_id_i}( \xbs_{i,-i}^k+\tau_i\textstyle\sum_{j=1 }^{N}w_{i,j}\xbs_{j,-i}^{k})
		\\
		\nonumber
		\0_{n_i}                             & \in  \partial_{x_i^{k+1}}  \bigl(
		\begin{aligned}[t]
		J_i(x_i^{k+1},\xbs_{i,-i}^{k+1}) +\textstyle \frac{1}{2\alpha\tau_i}\bigl\| x_i^{k+1}-x_i^{k} \bigr\|^2 \\
		+ \textstyle\frac{1}{2\alpha d_i}\bigl\| d_i x_i^{k+1}- \textstyle \sum_{j=1}^{N}w_{i,j}\xbs_{j,i}^k\bigr\|^2 \\
		+\iota_{\Omega_i}(x_i^{k+1})+\textstyle\frac{1}{\alpha}{(A_i^\top \lambda_i^k)}^\top x_i^{k+1}\bigr).
		\end{aligned}
		\end{aligned}
		\right.
		\end{align}
		Therefore, since the zeros of the  subdifferential of a (strongly) convex function coincide with the minima (unique minimum) \cite[Th.~16.3]{Bauschke2017}, \eqref{eq:useful1} can be rewritten as
		\begin{align}\label{eq:algorithmextform}
			\begin{aligned}
				&	\forall i\in \mc{I}:\! \left\{
					\begin{aligned}
						\xbs_{i,-i}^{k+1}    & = \textstyle \frac{1}{1+\tau_id_i}( \xbs_{i,-i}^k+\tau_i\textstyle\sum_{j=1 }^{N}w_{i,j}\xbs_{j,-i}^{k})
						\\
						x_i^{k+1}          & =
							\begin{aligned}[t]
								\underset{y \in \Omega_i} {\argmin}\bigl(
								J_i(y,\xbs_{i,-i}^{k+1})+\textstyle \frac{1}{2\alpha\tau_i}  \bigl\| y-x_i^{k} \bigr\|^2
								\\
								+\textstyle \frac{1}{2\alpha d_i}
								\bigl\| \textstyle d_i y-\sum_{j=1}^{N}w_{i,j}\xbs_{j,i}^{k} \bigr\|^2\\ +\textstyle\frac{1}{\alpha}{(A_i^\top \lambda_i^k)}^\top y  \bigr)
							\end{aligned}
					\end{aligned}
					\right.
					\\
					&
					\begin{aligned}
					\bs{v}^{k+1}       & = \bs{v}^{k}+{\bar\nu}\bs{V}_{\!\!m}\bs{\lambda}^k
					\\
					\bs{\lambda}^{k+1} & = \proj_{\R^{mN}_{\geq  0}}
					\bigl(
					\begin{multlined}[t]
						\bs{\lambda}^k +\bar{\delta}
						\bigl( \bs{A}\mc{R}(2\xbs^{k+1}-\xbs^k)-\bs{b}
						\\
						 ~~~~~~~~~~~~~~~~~~ -\bs{V}_{\!\!m}^\top (2\bs{v}^{k+1}-\bs{v}^k)
						\bigr)\bigr).
					\end{multlined}
					\end{aligned}
				\end{aligned}
		\end{align}
		The conclusion follows by defining $\bs{z}^k\coloneqq  \bs{V}_{\!\!m}^\top \bs{v}^k$, where $\bs{z}^k=\col((z_i)_{i\in\mc{I}})\in\R^{Nm}$ and $z_i^k\in \R^m$ are local variables kept by each agent, provided that $\bs{z}^0=\bs{V}_{\!\!m}^\top\bs{v}^0$. The latter is ensured by $\bs{z}^0=\0_{Nm}$, as in Algorithm~\ref{algo:1}.  \hfill $\blacksquare$
	\end{proof1}
	\begin{rem1}
		The preconditioning matrix $\Phi$ is
		designed to make the system in \eqref{eq:useful1} block triangular, i.e., to remove the term $\bs{W}_{\!\!n} \bs x^{k+1}$ and $\mc{R}^\top\bs{A}^\top \bs{\lambda}^{k+1}$ from the first inclusion, and the terms $\bs{V}_{\!\!m}\bs{\lambda}^{k+1}$ from the second one: in this way,    $\xbs_i^{k+1}$ and $\bs{z}^{k+1}$ do not depend on $\xbs_j^{k+1}$, for $i\neq j$, or $\bs{\lambda}^{k+1}$.
		This ensures that the resulting iteration can be computed by the agents  in a fully-distributed fashion (differently from the non-preconditioned resolvent $\J_{\mc{A}}$).
		Furthermore, the change of variable $\bs{z} = \bs{V}_{\!\! m}^\top \bs{v}$  reduces the number of auxiliary variables and decouples the dual update in \eqref{eq:algorithmextform} from the graph structure.
		\hfill $\square$
	\end{rem1}
	\begin{rem1}\label{rem:Singlevaluedfulldomain}
		By Lemma~\ref{lem:derivation}, Remark~\ref{rem:1_singlevaluedargmin}  and by  $\J_{\Phi^{-1}\Ac}$ in \eqref{eq:algorithmextform}, we conclude that  $\dom(\J_{\Phi^{-1}\Ac})=\R^{Nn+Em+Nm}$ and that $\J_{\Phi^{-1}\Ac}$ is single-valued.\hfill $\square$
	\end{rem1}

In order to apply Theorem~\ref{th:FNEconvergence} to the iteration in \eqref{eq:algcompact}, we still need the following lemma.
\begin{lem1}\label{lem:strmoninPhi}
   	Let $\alpha\in(0,\alphamax]$, $\alphamax$  as in Lemma~\ref{lem:strongmon_constant}. Then $\Phi^{-1} \Ac$ is restricted  monotone in $\mc{H}_\Phi$. \hfill $\square$.
\end{lem1}
\begin{proof1}
	Let $(\obs,\bs{u})\in\gra(\Phi^{-1}\mc{A})$, $\obs^*\in\zer(\Phi^{-1}\mc{A})$. Then, $(\obs,\Phi\bs{u}) \in \gra(\Ac)$ and $  \obs^*\in\zer(\mc{A})$. By  Lemma~\ref{lem:strongmon_constant} we conclude  that
	$ \langle \bs{u} \mid \obs-\obs^* \rangle_{\Phi} =\langle\Phi\bs{u}\mid \obs-\obs^* \rangle \geq 0.
	$
	\hfill$\blacksquare$
\end{proof1}

	\begin{thm}\label{th:main1}
		Let $\alpha\in(0,\alphamax]$, with $\alphamax$ as in Lemma~\ref{lem:strongmon_constant}, and let the step sizes $\bar{\tau},\bar{\nu},\bar{\delta}$ be as in Lemma~\ref{lem:stepsizes}. Then,  the sequence $(\xbs^k,\bs{z}^k,\bs{\lambda}^k)_{k\in\N}$ generated by Algorithm~\ref{algo:1} converges to some equilibrium $(\xbs^*,\bs{z}^*,\bs{\lambda}^*)$, where $\xbs^*=\1_N\otimes x^*$ and $x^*$ is the \gls{v-GNE} of the game in \eqref{eq:game}. \hfill $\square$
	\end{thm}
    \begin{proof1}
		By Lemma~\ref{lem:derivation}, we can equivalently study the convergence of the iteration in \eqref{eq:algcompact}. In turn, \eqref{eq:algcompact} can be rewritten as \eqref{eq:KM} with  $\gamma^k=1$, $e^k=\0$, for all $k\in\N$. For later reference, let us define $\ubs^k=\J_{\Phi^{-1} \mc{A}}(\obs^k)$ (here $\ubs^k=\obs^{k+1}$).  ${\Phi^{-1} \mc{A}}$ is restricted monotone in $\mc{H}_{\Phi}$ by Lemma~\ref{lem:strmoninPhi}.
		By Theorem~\ref{th:FNEconvergencei}, the sequence $(\obs^k)_{k\in\N}$ is bounded, hence it admits at least one cluster point, say $\bar{\obs}$.  By \eqref{eq:explicitinclusion} and \eqref{eq:opA}, it holds, for any $\obs \in \bs{\Omega}\times \R^{Em} \times \R_{\geq 0}^{Nm}$, that $\langle \Ac_1(\ubs^{k})+\Phi(\ubs^{k}-\obs^k) \mid \obs-\ubs^{k}\rangle \geq 0 $, with $\Ac_1$ as in \eqref{eq:opA}. By Theorem~\ref{th:FNEconvergenceii}, $\ubs^{k}-\obs^k \rightarrow \0$. Therefore, by continuity of $\Ac_1$, taking the limit on a  diverging  subsequence $(l_k)_{\k}$ such that $(\obs^{l_k})_\k \rightarrow \obsbar$, we have that for all
		$
 		 \obs \in \bs{\Omega}\times \R^{Em} \times \R_{\geq 0}^{Nm}),$ 	$\langle \Ac_1(\obsbar)\mid \obs-\obsbar\rangle \geq 0,
		$
		 which shows that $\obsbar\in \zer(\Ac)=\fix(\J_{\Phi^{-1} \mc{A}})$.
		 Hence
		$(\obs^k)_\k$ converges to an equilibrium of \eqref{eq:algcompact}  by  Theorem~\ref{th:FNEconvergenceiii}. The conclusion follows by Lemma~\ref{lem:zeros}. \hfill  $\blacksquare$
	\end{proof1}
\begin{rem1}
		While the choice of  step sizes in Lemma~\ref{lem:stepsizes} is decentralized, computing the bound $\alphamax$ for the common parameter $\alpha$ in Algorithm~\ref{algo:1} requires some global information on the graph $\mc{G}$ (i.e., the algebraic connectivity) and on the game mapping (the strong monotonicity and Lipschitz constants).
\hfill $\square$
\end{rem1}
\begin{rem1}
		If $\xbs^0\in\Omega^N$, then $\xbs^k\in \Omega^N$ for all $\k$ (by convexity and the updates in Algorithm~\ref{algo:1}), and  Assumption~\ref{Ass:StrMon} can be relaxed to hold only on $\Omega$.  \hfill $\square$
\end{rem1}
	\begin{rem1}[{\it Inexact updates}] \label{rem:reminexact}
		 The local optimization problems in Algorithm~\ref{algo:1} are strongly convex, hence they can be efficiently solved by several iterative algorithms (with linear rate). While computing the exact solutions $\bar{x}_i^k$ would require an infinite number of iterations, the convergence in Theorem~\ref{th:main1} still holds if $x_i$ is updated with an approximation $\widehat{x}_i^k$ of $\bar{x}_i^k$, provided that the errors  $e^k_i\coloneqq \bar{x}^k_i-\widehat{x}_i^k$ are norm summable, i.e.,  $(\|e_i^k\|)_{k\in\N}\in\ell^1$, for all $i \in \mc{I}$ (the same proof applies, since the condition on $e^k$  in Theorem~\ref{th:FNEconvergence} would be satisfied, by equivalence of norms). For example, assume that $\widehat{x}_i^k$ is computed via a finite number $j_i^k\geq 1$ of steps of the projected gradient method, warm-started at $x_i^k$, with (small enough) fixed step.
		 Then, each  agent  can independently ensure that $\|e_i^k\|\leq \varepsilon_i^k$,  for some $(\varepsilon_i^k)_{\k}\in\ell^1$,
		 by simply choosing
		 \begin{equation}\label{eq:stopping criterion}
		 j_i^k \geq \log \left (  \varepsilon_i^k(1-\rho_i)/{\|x_i^k-\widehat{x}^{k,1}_i\|} \right) / \log(\rho_i),
		 \end{equation}
		 where $\widehat{x}^{k,1}_i$ is the approximation obtained after one gradient step and $\rho_i\in(0,1)$ is the contractivity parameter of the gradient descent\footnote{{$\rho_i$ can be taken independent of $k$: since $\nabla J_i(\cdot,\xbs_{i,-i})$ is $\mu_i$ strongly monotone and $\theta_i$ Lipschitz, for some $\mu_i\geq\mu$, $\theta_i \leq \theta$ and  for all  $\xbs_{i,-i}$,  the factor $\rho_i= \frac{\theta_i -\mu_i}{\theta_i +\mu_i+1/(\alpha \tau_i) + d_i/\alpha}$ is ensured by the step  $2/({\theta_i +\mu_i+1/(\alpha \tau_i) + d_i/\alpha})$}.}.
		 We finally remark that $\bar{x}_i^k$ must be estimated with increasing accuracy.  In practice, however, when $x_i^k$ is converging, $\| x_i^{k+1}-x_i^{k}\|\rightarrow 0$. Hence $x_i^k$ is a good initial guess for $\bar{x}_i^k$, and the computation of $x_i^{k+1}$ often requires few gradient steps, see also §\ref{sec:numerics}.
		\hfill $\square$
	\end{rem1}
\section{Accelerations }\label{sec:acc}
\begin{algorithm*}\caption{Fully-distributed \gls{v-GNE} seeking via accelerated \gls{PPPA}  } \label{algo:2_overrel}
	\vspace{0.5em}
	\begin{itemize}[leftmargin=6.3em]
		\item[Initialization:]
		\begin{itemize}[leftmargin=1em]
			\item
			Choose acceleration:  \begin{tabular}{|ll }
				Overrelaxation:& set $\gamma>0$, $\zeta=0$, $\eta=0$; \\[-5pt]
				Inertia:  &   set $\gamma=0$, $\zeta>0$, $\eta=0$;    \\[-5pt]
				Alternated inertia: & set $\gamma=0$, $\zeta=0$, $\eta>0$; \end{tabular}
			\newline
			\item For all $i\in \mc{I}$, set $x_i^{-1}=x_i^0\in \Omega_i$, $\xbs_{i,-i}^{-1}=\bs{x}_{i,-i}^0\in \R^{n-n_i}$, $z_i^{-1}=z_i^0=\0_m$, $\lambda_i^{-1}=\lambda_{i}^0\in \R^m_{\geq 0}$.
		\end{itemize}
		\vspace{1em}
		\item[For all $k> 0$:]
		\begin{itemize}[leftmargin=1em]
			\item (Alternated) inertial step:  set $\tilde{\eta}^k =0$  if $k$ is even, $ \tilde{\eta} ^k=\eta$ otherwise; each agent $i\in\mc{I}$ computes
			\begin{align*}
				\vspace{0.4em}
				\tilde{\xbs}_{i,-i}^{k} &= {\xbs}_{i,-i}^{k} +(\zeta+\tilde{\eta}^k) ({\xbs}_{i,-i}^{k}-{\xbs}_{i,-i}^{k-1})
				&   &   &
				\tilde{x}_i^{k} &= {x}_{i}^{k} +(\zeta+\tilde{\eta}^k)({x}_{i}^{k}-{x}_{i}^{k-1})
				\\
				\tilde{z}_{i}^{k} &= {z}_{i}^{k} +(\zeta+\tilde{\eta}^k)({z}_{i}^{k}-{z}_{i}^{k-1})
				&   &   &
				\tilde{\lambda}_{i}^{k} &= {\lambda}_{i}^{k} +(\zeta+\tilde{\eta}^k)({\lambda}_{i}^{k}-{\lambda}_{i}^{k-1})
			\end{align*}
			\item
			Communication: The  agents exchange the variables  $\{ \tilde x^k_{i},\tilde{\bs{x}}_{i,-i}^k, \tilde \lambda_i^k \}$ with their neighbors.
			\item
			Resolvent computation:  each agent $i\in\mc{I}$ computes \vspace{-0.3em}
			\begin{align*}
				\vspace{0.4em}
				\breve{\xbs}_{i,-i}^{k+1} & = \textstyle \frac{1}{1+\tau_id_i}( \tilde\xbs_{i,-i}^k+\tau_i\textstyle\sum_{j=1 }^{N}w_{i,j}\tilde\xbs_{j,-i}^{k})
				\\
				\vspace{0.4em}
				\breve{x}_{i}^{k+1}     & =
				\underset{y \in \Omega_i} {\argmin}\bigl(
				J_i(y,\breve\xbs_{i,-i}^{k+1})+\textstyle \frac{1}{2\alpha\tau_i}  \bigl\|  y-\tilde x_i^{k} \bigr\|^2
				+\textstyle \frac{1}{2\alpha d_i}
				\bigl\| \textstyle d_i y-\sum_{j=1}^{N}w_{i,j}\tilde\xbs_{j,i}^{k} \bigr\|^2+\frac{1}{\alpha}(A_i^\top \tilde\lambda_i^k)^\top y  \bigr)
				\\
				\vspace{0.4em}
				\breve{z_i}^{k+1}       & = \tilde z_i^{k}+\textstyle \sum_{j=1}^{N} \nu_{(i,j)} w_{i,j} (\tilde\lambda_i^k-\tilde\lambda_j^k)
				\\
				\vspace{0.4em}
				\breve{\lambda}_i^{k+1} & =\proj_{\R^{m}_{\geq  0}}
				\bigl(
				{\tilde\lambda}_i^k +\textstyle\delta_i
				\bigl( A_i(2 \breve x_i ^{k+1}-\tilde x_i^k)-b_i-(2 \breve z_i^{k+1}-\tilde z_i^k)
				\bigr)
				\bigr).
			\end{align*}
			\item Relaxation step: each agent $i\in\mc{I}$ computes \vspace{-0.3em}
			\begin{align*}
				\xbs_{i,-i}^{k+1} &= \gamma  \breve{\xbs}_{i,-i}^{k+1} +(1-\gamma){\xbs}_{i,-i}^{k}
				&   &   &
				x_i^{k+1} & =\gamma  \breve{x}_{i}^{k+1} +(1-\gamma){x}_{i}^{k}
				\\
				{z_i}^{k+1} &=\gamma  \breve{z}_{i}^{k+1} +(1-\gamma){z}_{i}^{k}
				&   &   &
				{\lambda}_i^{k+1} & =\gamma  \breve{\lambda}_{i}^{k+1} +(1-\gamma){\lambda}_{i}^{k}
			\end{align*}
			\vspace{-1em}
		\end{itemize}
	\end{itemize}
\end{algorithm*}
	Lemma~\ref{lem:derivation} shows that Algorithm~\ref{algo:1} can be recast (modulo the  change of variables $\bs{z}=\bs{V}_{\!\!m}^\top\bs{v}$) as
	\begin{equation}
	\obs^{k+1}=T(\obs^k),
	\end{equation}
	where $T\coloneqq \J_{\Phi^{-1}\Ac}$. This compact operator representation allows for some modifications of Algorithm~\ref{algo:1}, that can increase its convergence speed. In particular, we consider three popular accelerations schemes \cite{Iutzeler_Hendrickx2019}, which have been extensively studied for the case of firmly nonexpansive operators \cite[Def.~4.1]{Bauschke2017}, and also found application in games under full-decision information \cite{Belgioioso_aggregative_2020}, \cite{Scutari_ComplexGames_TIT2014}. Here  we provide convergence guarantees for the  partial-decision information setup, where  $T$ is only firmly quasinonexpansive. Our fully distributed accelerated algorithms are illustrated in   Algorithm~\ref{algo:2_overrel}.  In the following, we  assume that   $\alpha\in(0,\alphamax]$,  $\alphamax$ as in Lemma~\ref{lem:strongmon_constant}, and that the step sizes $\bar{\tau},\bar{\nu},\bar{\delta}$ are chosen as in Lemma~\ref{lem:stepsizes}.
	\begin{prop1}[Overrelaxation] \label{prop:o}
	Let $\gamma \in [1,2)$. Then, for any $\obs^0$, the sequence $(\obs^k)_\k$ generated by
	\begin{equation}\label{eq:overrelax}
	(\forall \k), \quad  \obs^{k+1}=\obs^k+\gamma (T(\obs^k)-\obs^k),
	\end{equation}
	converges to an equilibrium $(\xbs^*,\bs{v}^*,\bs{\lambda}^*)\in\zer(\Ac)$, where $\xbs^*=\1_N\otimes x^*$ and $x^*$ is the \gls{v-GNE} of the game in \eqref{eq:game}.
	\hfill $\square$
	\end{prop1}

	\begin{proof1}
		The iteration in \eqref{eq:overrelax} is in the form \eqref{eq:KM}, with $\gamma^k=\gamma$, $e^k=\0$, for all $\k$. Then, the conclusion follows analogously to Theorem~\ref{th:main1}. \hfill $\blacksquare$
	\end{proof1}
	\begin{prop1}[Inertia] \label{prop:i}
		Let $\zeta \in [0,\frac{1}{3})$. Then, for any $\obs^{-1}\coloneqq \obs^0$, the sequence $(\obs^k)_\k$ generated by
		\begin{equation}\label{eq:inertia}
		(\forall \k), \quad  \obs^{k+1}=T(\obs^k+\zeta(\obs^k-\obs^{k-1})),
		\end{equation}
		converges to an equilibrium $(\xbs^*,\bs{v}^*,\bs{\lambda}^*)\in\zer(\Ac)$, where $\xbs^*=\1_N\otimes x^*$ and $x^*$ is the \gls{v-GNE} of the game in \eqref{eq:game}.
		\hfill $\square$
	\end{prop1}

	\begin{proof1}[Proof (sketch)]
		By following all the steps in the proof of \cite[Th.~5]{Bot2015_InertialDR} (which can be done by recalling that  an operator $T$ is firmly (quasi)nonexpansive if and only if the operator $2T-\Id$ is (quasi)nonexpansive \cite[Prop. 4.2, 4.4]{Bauschke2017}),
it can be shown that, if $\zeta\in[0,\frac{1}{3})$, then $(\obs^k)_\k$ is bounded and $\obs^{k+1}-\obs^{k}\rightarrow 0$. Then, the proof follows analogously to Theorem~\ref{th:main1}. \hfill $\blacksquare$
	\end{proof1}

		\begin{prop1}[{Alternated inertia}] \label{prop:ai}
		Let $\eta \in [0,1]$. Then, for any $\obs^0$, the sequence $(\obs^k)_\k$ generated by
		\begin{equation} \label{eq:alternated_inertia}
		\left\{
		\begin{aligned}
		\obs^{k+1} & =T (\obs^k)                     & \text{if $k$ is even,} \\
		\obs^{k+1} & =T (\obs^k+\eta(\obs^k-\obs^{k-1})) & \text{if $k$ is odd,}
		\end{aligned}
		\right. \hspace{-0.1em}
		\end{equation}
		converges to an equilibrium $(\xbs^*,\bs{v}^*,\bs{\lambda}^*)\in\zer(\Ac)$, where $\xbs^*=\1_N\otimes x^*$ and $x^*$ is the \gls{v-GNE} of the game in \eqref{eq:game}.
		\hfill $\square$
	\end{prop1}
   \begin{proof1}
   	For all $\k$, $\obs^{2k+2}=T(T(\obs^{2k})+\eta(T(\obs^{2k})-\obs^{2k})),$
   	which is the same two-steps update obtained in  \eqref{eq:KM} with $\gamma^{2k}=1+\eta$, $\gamma^{2k+1}=1$ (and $\mc{B}\coloneqq \Phi^{-1}\Ac$, $e^k=\0$). Therefore the convergence of the sequence $(\obs^{2k})_\k$ to an equilibrium $(\xbs^*,\bs{v}^*,\bs{\lambda}^*)\in\zer(\Ac)$ follows analogously to Theorem~\ref{th:main1} (with a minor modification for the case $\eta=1$). The convergence of the sequence $(\obs^{2k+1})_\k$ then follows by Theorem~\ref{th:FNEconvergencei}.
   	\hfill $\blacksquare$
   	\end{proof1}

We note that, by Theorem~\ref{th:FNEconvergence}, the convergence results in Propositions~\ref{prop:o} and \ref{prop:ai} hold also in the case of summable errors on the updates, as in Remark~\ref{rem:reminexact}. Analogously to our analysis,  provably convergent acceleration schemes could also be obtained for the \gls{FB} algorithm in \cite{Pavel2018}:  however, an advantage of our \gls{PPA} is  that the bounds on the inertial/relaxation parameters are fixed and independent on (unknown) problem parameters.

\subsection{On the convergence rate}\label{subsec:rate}
We conclude this section with a discussion on the convergence rate of Algorithms~\ref{algo:1} and\ref{algo:2_overrel}.
First, even  under Standing Assumption~\ref{Ass:StrMon}, the \gls{KKT} operator on the right-hand side of \eqref{eq:KKT} is generally not strongly monotone. Similarly, the operator $\mc{A}$ in \eqref{eq:opA} is not strongly monotone and Algorithm~\ref{algo:1} can have multiple fixed points.
Therefore, one  should not expect linear convergence.
By Lemma~\ref{lem:derivation} and the proof of Theorem~\ref{th:FNEconvergence}, we can derive the following ergodic rate for the fixed-point residual in Algorithm~\ref{algo:1}:
\[
\textstyle
\frac{1}{k} \sum_{i=0}^{k}  \| \obs^{k+1}-\obs^k\|^{2} \leq {O(1/k).}
\]
This rate also holds for the iterations in \eqref{eq:overrelax}, \eqref{eq:inertia}, \eqref{eq:alternated_inertia}; for the case of general operator splittings (and differently from optimization algorithms), tighter  rates for accelerated schemes are only known for particular cases, and most works focus on mere convergence \cite{Iutzeler_Hendrickx2019}, \cite{Bot2015_InertialDR}. Yet, the practice shows that relaxation and inertia often result in improved speed, see \cite{Belgioioso_aggregative_2020} or §\ref{sec:numerics}.

The same residual rate $O(1/k)$ can also be shown for the  pseudo-gradient method in \cite[Alg.~1]{Pavel2018}. However, a major difference from Lemma~\ref{lem:stepsizes} is that the upper bounds for the step sizes in \cite[Th.~2]{Pavel2018} are proportional to the constant ${\mu_{\Fa}}$  in \eqref{eq:alphamax}, which is typically very small (up to scaling of the whole operator $\Fa$),
\cite{Bianchi_LCSS2020} (see also §\ref{subsec:NC}), and, most importantly, it vanishes as the number of agents increases (fixed the other parameters).
In contrast, our algorithms allows for much larger  steps,  which can be chosen independently of the number of agents. {This is a structural advantage of the \gls{PPA}, whose convergence does not depend on the cocoercivity constant of the operators involved.}
Indeed, step sizes must  be taken into account if  convergence is evaluated in terms of residuals.

We finally note that linear convergence can be achieved via \gls{PPPA} for games without coupling constraints.
For instance, Algorithm~\ref{algo:1NE} corresponds to the overrelaxed method in Algorithm~\ref{algo:2_overrel}, and can be derived, as in Lemma~\ref{lem:derivation}, by taking $\mc{B}=\Phi_{\textnormal{NE}}^{-1}\mc{A}_{\textnormal{NE}}(\bs{x})$ in \eqref{eq:KM}, where
$\mc{A}_{\textnormal{NE}}(\bs{x})\coloneqq \Fa(\bs{x})+\nc_{\bs{{\Omega}}}(\bs{x})$ and  $\Phi_{\textnormal{NE}}\coloneqq \bar{\tau}^{-1}+\bs{W}_{\!\!n}$ are obtained by removing the dual variables from $\mc{A}$, $\Phi$.
By   \eqref{eq:strongFa}, as in Lemma~\ref{lem:strmoninPhi}, it can be shown that $\mc{A}_{\textnormal{NE}}$ is restricted $\frac{\mu_{\Fa} }{\|\Phi_{\textnormal{NE}}\|}$-strongly monotone in $\mc{H}_{\Phi_{\textnormal{NE}}}$. Thus, recursively applying Theorem~\ref{th:FNEconvergenceiv}, we can infer the following result, which appeared in \cite{Bianchi_CDC20_PPP} only limited to  $\gamma=1$.
\begin{algorithm}[t] \caption{Fully-distributed NE seeking via  \gls{PPPA} } \label{algo:1NE}
	\vspace{-0.4em}
	\begin{align*}
		\breve{\xbs}_{i,-i}^{k+1} & = \textstyle \frac{1}{1+\tau_id_i}( \xbs_{i,-i}^k+\tau_i\textstyle\sum_{j=1 }^{N}w_{i,j}\xbs_{j,-i}^{k})
		\\
		\breve{x}_{i}^{k+1}     & =
		\underset{y \in \Omega_i} {\argmin}\bigl(
		J_i(y,{\breve\xbs}_{i,-i}^{k+1})+\textstyle \frac{1}{2\alpha\tau_i}  \bigl\|  y- x_i^{k} \bigr\|^2
		\\[-0.5em]
		& \qquad \qquad \qquad +\textstyle \frac{1}{2\alpha d_i}
		\bigl\| \textstyle d_i y-\sum_{j=1}^{N}w_{i,j} \xbs_{j,i}^{k} \bigr\|^2 \bigr)
		\\[0.5em]
		\xbs_i^{k+1}& =\xbs_i^{k}+\gamma(\breve{\xbs}_{i}^{k+1} - \xbs_i^{k})
	\end{align*}
	\vspace{-1em}
\end{algorithm}
\begin{thm}
	\label{th:mainNE}
	Let $\tau_i^{-1}>d_i$ for all $i\in\mc{I}$, let $\gamma\in(0,2)$, and let $\alpha\in(0,\alphamax]$, with $\alphamax$ as in Lemma~\ref{lem:strongmon_constant}.
	Then, the sequence $(\xbs^k)_{k\in\N}$ generated by Algorithm~\ref{algo:1NE} converges to $\xbs^*=\1_{N}\otimes x^*$, where $x^*$ is the unique Nash equilibrium of the game in \eqref{eq:game}, with linear rate:
	\[ (\forall \k) \quad \|\xbs^k-\xbs^*\|_{\Phi_{\textnormal{NE}}}\leq \textstyle (\rho_\gamma)^{ k}\|\xbs^0-\xbs^*\|_{\Phi_{\textnormal{NE}}}, \]
	where $\rho_\gamma\coloneqq \max(1-\textstyle \frac{\gamma\mu_{\Fa} }{\|\Phi_{\textnormal{NE}}\|+\mu_{\Fa}} ,\gamma-1)$, $\mu_{\Fa}$ as in \eqref{eq:alphamax}. \hfill $\square$
	%
\end{thm}
The best theoretical rate
$\rho_{\bar{\gamma}}=1-{2\mu_{\Fa}}/({\|\Phi_{\textnormal{NE}}\|+2\mu_{\Fa}})$
is obtained for  $\bar{\gamma}=1+{\|\Phi_{\textnormal{NE}}\|}/({\|\Phi_{\textnormal{NE}}\|+2\mu_{\Fa}})$.
We  observed in \cite[§5]{Bianchi_CDC20_PPP} that this rate compares favorably with that of the  state-of-the-art algorithms -- please refer to \cite{Bianchi_CDC20_PPP}, also for numerical results. For instance, in the absence of coupling constraints, the \gls{FB} algorithm in \cite[Alg.~1]{Pavel2018} reduces to \cite[Alg.~1]{TatarenkoShiNedicTAC20}, whose optimal linear rate  $O ( (1-{\kappa_{\Fa}}^2)^\frac{k}{2} )$ depends \emph{quadratically} on the quantity $\kappa_{\Fa}\coloneqq\mu_{\Fa} / \theta_{\Fa} < 1$ \cite[Th.~7]{TatarenkoShiNedicTAC20}, where $\theta_{\Fa}\coloneqq 2\max((d_i)_{i\in\mc{I}})+\alpha\theta$. Instead,  $\rho_{\bar{\gamma}}\leq1-\kappa_{\Fa}$, for large enough  $\tau_i$'s (since ${\|\Phi_{\textnormal{NE}}\|+2\mu_{\Fa}}\leq \max((d_i+\tau_i^{-1})_{i\in\mc{I}})+2\alpha\theta$), as shown in Table~\ref{tab:1}.
\begin{table}[t]
	\centering
	\newlength{\mywidth}
	\setlength{\mywidth}{0.8em}
	\begin{tabular}{  c c c  }
		\toprule
		&
		\gls{FB} \cite[Alg.~1]{Pavel2018}
		&
		\hspace{0.2em} \gls{PPPA} \hspace{0.2em}
		\\
		\midrule
		step sizes
		&
		\begin{tabular}{c}  $O\left(\frac{\mu_{\Fa}}{{\theta_{\Fa}}^2+\mu_{\Fa}} \right)$ \end{tabular}
		&
		$O(1)$
		\\
		\begin{tabular}{c} linear rate $\rho$ (no \\[-0.5em]  coupling constraints)  \end{tabular}
		&
		$(1-{\kappa_{\Fa}}^2) ^{\frac{1}{2} }$
		&
		$ 1-{\kappa_{\Fa}} $
		\\  \bottomrule
	\end{tabular}
	\vspace{0.3em}
	\caption{\label{tab:1}Comparison between our \gls{PPPA} and  projected pseudo-gradient methods.}
\end{table}
		\section{Aggregative games}\label{sec:aggregative}
			\begin{algorithm*} \caption{Fully-distributed \gls{v-GNE} seeking in aggregative games via \gls{PPPA} } \label{algo:1agg}
			\vspace{0.4em}
			\begin{itemize}[leftmargin=6.3em]
				\item[Initialization:]
				\begin{itemize}[leftmargin=1em]
					\item For all $i\in \mc{I}$, set $x_i^0\in \Omega_i$, $s_i^0=\0_{\bar{n}}$, $z_i^0=\0_m$, $\lambda_{i}^0\in \R^m_{\geq 0}$.
				\end{itemize}
				\vspace{0.1em}
				\item[For all $k> 0$:]
				\begin{itemize}[leftmargin=1em]
					\item
					Communication: The  agents exchange the variables  $\{ \sigma_i^k=x^k_{i}+s_i^k,\lambda_i^k \}$ with their neighbors.
					\item
					Local variables update: each agent $i\in\mc{I}$ computes
					\begin{align}
					\nonumber
					\vspace{0.4em}
					s_i^{k+1}         & =s_i^k -\textstyle \beta \sum_{j=1}^N w_{i,j}({\sigma_i^k- \sigma_j^k})
					\\
					\label{eq:xupdateagg}
					\vspace{0.4em}
					x_i^{k+1}           & \leftarrow y \ \   \text{s.t.}  \ \  \0_{\bar{n}} \in \alpha \bs{\tilde{F}}_{\! i}(y,y+s_i^{k+1})
					+\textstyle \frac{1}{\tau_i} (y-x_i^k)+A_i^\top {\lambda}_i^k+ \textstyle \sum_{j=1}^{N}w_{i,j}(\sigma_i^k-\sigma_j^k) +\mathrm{N}_{{\Omega_i}}(y) \!\!
					%
					%
					%
					\\
					\nonumber
					\vspace{0.4em}
					{z_i}^{k+1}       & = z_i^{k}+\textstyle \sum_{j=1}^{N} \nu_{(i,j)} w_{i,j} (\lambda_i^k-\lambda_j^k)
					\\
					\nonumber
					\vspace{0.4em}
					{\lambda}_i^{k+1} & =\proj_{\R^{m}_{\geq  0}}
					\bigl(
					{\lambda}_i^k +\textstyle \delta_i
					\bigl( A_i(2x_i^{k+1}-x_i^k)-b_i-(2z_i^{k+1}-z_i^k)
					\bigr)
					\bigr).
					\end{align}
				\end{itemize}
			\end{itemize}
			\vspace{-0.6em}
		\end{algorithm*}
	In this section we focus on the  particularly relevant class of  (average) aggregative games, which arises in a variety of engineering applications, e.g.,  network congestion control and demand-side management \cite{Grammatico2017}. In aggregative games,  $n_i=\bar{n}>0 $ for all $i\in\mc{I}$ (hence $n=N\bar{n}$) and the cost function of each agent depends only on its local decision and on the value of the average strategy
	$
	{\operatorname{avg}}(x)\coloneqq  \tfrac{1}{N}\textstyle{ \sum_{i\in\mc{I}} }x_i.
	$
	Therefore, for each $i\in \mc{I}$, there is a function $f_i:\R^{\bar{n}}\times \R^{\bar{n}}\rightarrow {\R}$ such that the original cost function $J_i$ in \eqref{eq:game} can be written as
	\begin{align}\label{eq:Jf_aggregative}
	J_i(x_i,x_{-i})=:f_i(x_i,\operatorname{avg}(x)).
	\end{align}
	Since an aggregative game is only a particular instance of the game in \eqref{eq:game}, all the considerations
	on  the existence and  uniqueness of a  \gls{v-GNE} and the equivalence with the KKT conditions in \eqref{eq:KKT}  are still valid.

	Moreover, Algorithms~\ref{algo:1} could still be used to compute a \gls{v-GNE}. This would require each agent to keep (and exchange) an estimate of all other agents' action, i.e., a vector of $(N-1)\bar{n}$ components. In practice, however, the cost of each agent is only a function of the aggregative value $\operatorname{avg}(x)$, whose dimension $\bar{n}$ is independent of the number $N$ of agents. To reduce  communication and computation burden,  in this section we introduce a \gls{PPPA} specifically tailored to seek a \gls{v-GNE} in aggregative games, that is scalable with the number of agents. The proposed iteration is shown in Algorithm~\ref{algo:1agg}, where the parameters $\alpha$, $\beta$, and  $\tau_i$, $\delta_i$ for all $i \in \mc{I}$, $v_{(i,j)} $ for all $(i,j)\in\mc{E}$  have to be chosen appropriately, and  we denote
	\begin{align}
	\label{eq:Fitilde}
	\Fi(x_i,\xi_i)\coloneqq  \nabla_{\!x_{i}} f_i(x_{i}, \xi_i)+ \textstyle \frac{1}{N}\nabla_{\! \xi_i}f_i(x_{i}, \xi_{i}).
	\end{align}
	We note that $\Fi(x_i, \avg(x)) =\nabla_{\!x_i}J_i(x_i, x_{-i})= \nabla_{\!x_i}f_i(x_i, \avg(x))$.

		Because of the partial-decision information assumption, no agent has access to  the actual value of the average strategy. Instead, we equip each agent with an auxiliary error variable $s_i\in \R^{\bar{n}}$, which is an estimate of the quantity $\operatorname{avg}(x)-x_i$.
		Each agent aims at reconstructing the true aggregate value, based on the information received from its neighbors.
	In particular, it should hold that $s^k\rightarrow \1_N\otimes \avg(x^k)-x^k$ asymptotically, where $s\coloneqq \col((s_i)_{i\in\mc{I}})$. For brevity of notation, we  also denote
	\begin{equation}\label{eq:sigmadef}
	\sigma_i\coloneqq x_i+s_i, \quad  \sigma\coloneqq \col((\sigma_i)_{i\in\mc{I}}).
	\end{equation}
       \begin{rem1}\label{rem:invariance}
 		By the updates in Algorithm~\ref{algo:1agg}, we can infer an important invariance property, namely that $\avg(s^k)=\0_{\bar{n}}$, or equivalently  $\avg(x^k)=\avg(\sigma^k)$, for any $k\in\N$, provided that the algorithm is initialized appropriately, i.e., $s_i^0=\0_{\bar{n}}$, for all $i\in \mc{I}$.
		 In fact, the update of $\sigma$, as it follows from Algorithm~\ref{algo:1agg}, is
		 \begin{equation}\label{eq:dyn_track}
		 \sigma^{k+1} =\sigma^k-\textstyle \beta\bs{L}_{\bar{n}}\sigma^k +(x^{k+1}-x^k),
		 \end{equation}
		where   $\bs{L}_{\bar{n}}\coloneqq L\otimes I_{\bar{n}}$.
		This update  is a dynamic tracking for the time-varying quantity $\avg(x)$, similar to those considered for aggregative games in \cite{Koshal_Nedic_Shanbag_2016}, \cite{BelgioiosoNedicGrammatico2020}, \cite{GadjovPavel_aggregative_2019}. Differently from \cite{GadjovPavel_aggregative_2019}, here we introduce the  error variables $s_i$, which allow us to directly recast the iteration in \eqref{eq:dyn_track} in an operator-theoretic framework.
	    \hfill $\square$
\end{rem1}
Similarly to §\ref{sec:derivationconvergence}, we study the convergence of Algorithm~\ref{algo:1agg} by relating it to the iteration in \eqref{eq:KM}.
First, let us define the \emph{extended pseudo-gradient} mapping
\begin{align}
\label{eq:extended_pseudo-gradient_agg}
\bs{\tilde{F}}(x,\xi)\coloneqq
\col((\bs{\tilde{F}}_{\! i}(x_i,\xi_i))_{i\in\mc{I}} ),
\end{align}
with  $\xi\coloneqq \col((\xi_i)_{i\in\mc{I}})\in \R^n$, and the operators  $\tilde{\Fa}(x,s)\coloneqq \col(\alpha\bs{\tilde{{F}}}(x, {\sigma})+\bs{L}_{\bar{n}} {\sigma}, \bs{L}_{\bar{n}}{\sigma} ) $,
\begin{align}
\label{eq:opAtilde}
\tilde{\mathcal{A}}(\obs) \! & \coloneqq
 \! \hspace{-0.2em}
\left[
\begin{matrix}
\alpha\bs{\tilde{{F}}}(x, {\sigma})+\bs{L}_{\bar{n}} {\sigma}\\
\bs{L}_{\bar{n}}{\sigma}\\
\0_{Em}\\
\bs b
\end{matrix}\right] \!
\hspace{-0.2em}+ \hspace{-0.2em}
\left[
\begin{matrix}
{\bs{A}^\top \bs{\lambda}} \\ \0_n \\ {-\bs{V}_{\!\! m} \bs{\lambda}} \\ {{\bs{V}_{\!\! m}^\top \bs{v} -\bs{A}x  }}
\end{matrix}\right]
\! \hspace{-0.2em} + \! \hspace{-0.2em}
\left[
\begin{matrix}
\nc_{{\Omega}}(x) \\
\0_{n}\\
\0_{Em}\\
\nc_{\R^{Nm}_{\geq 0}}(\bs{\lambda})
\end{matrix}\right]\!\!,
\end{align}
where
$\obs\coloneqq \col(x,s,\bs{v},\bs{\lambda})\in\R^{2n+Em+Nm}$,
and we recall  that ${\sigma}=x+s$ is just a shorthand notation.
	\begin{lem1}\label{lem:LipschitzExtPseudoagg}
	The   mapping $\bs{{\tilde F}}$ in \eqref{eq:extended_pseudo-gradient_agg} is $\tilde{\theta}$-Lipschitz continuous, for some $\tilde{\theta}>0$.
	{\hfill $\square$} \end{lem1}
\begin{proof1}
	It follows from Lemma~\ref{lem:LipschitzExtPseudo}, by noticing that $\bs{\tilde{F}}(x,\xi)=\bs{F}((x,(I_N\otimes\1_{N-1}\otimes I_{\bar{n}})(\textstyle \frac{N}{N-1}{\xi}-\frac{1}{N-1}x)))$. \hfill $\blacksquare$
\end{proof1}
Finally, we will assume that the  step sizes $\bar{\tau}\coloneqq\diag((\tau_i I_{\bar{n}})_{i\in\mc{I}})$, $\bar{\beta}\coloneqq\beta I_{Nn}$, $\bar{\nu}\coloneqq\diag((\nu_{(i,j)})_{(i,j)\in\mc{E}})$, $\bar{\delta}\coloneqq\diag((\delta_i I_ m)_{i\in\mc{I}})$ are  chosen such that  $\tilde\Phi \succ 0$, where
\begin{equation}\label{eq:Phitilde}
\tilde \Phi\coloneqq  \left[
\begin{matrix}
{\bar{\tau}^{-1}-\bs{L}_{\bar{n}}} & -\bs{L}_{\bar{n}}            & {\0}                     & -{\bs{A}^\top}
\\
-\bs{L}_{\bar{n}}             & \bar{\beta}^{-1}-\bs{L}_{\bar{n}} & \0                       & \0
\\
\0                            & \0                           & {\bar{\nu}}^{-1}            & {\bs{V}_{\!\! m} }
\\
{-\bs{A}}                     & \0                           & {\bs{V}_{\!\! m} ^\top } & {\bar{\delta}}^{-1}
\end{matrix}\right].
\end{equation}
\begin{lem1}\label{lem:stepsizesagg}
	The matrix $\tilde{\Phi}$ in \eqref{eq:Phitilde} is positive definite if  $\beta^{-1}  > 4\max((d_i)_{i\in\mc{I}})$, ${\nu_{(i,j)}}^{-1} >  2\sqrt{(w_{i,j})}$ for all $(i,j)\in\mc{E}$, and
	$\tau_i ^{-1}> 4d_{i} +\|A_i^\top\|_{\infty}$,
	$\delta_i ^{-1} > \|A_i \|_{\infty}+\textstyle \sum_{j=1}^N \sqrt{w_{i,j}}$ for all $i\in \mc{I}$. \hfill $\square$
\end{lem1}
\begin{thm}\label{th:main1agg}
	Let $d_{\min}\coloneqq \min((d_i)_{i\in{\mc{I}}})$ and
	\begin{equation}\label{eq:alphamaxagg}
	{\tilde{\alpha}_{\textnormal{max}}}\coloneqq \min\left( \textstyle \frac{4\mu \uplambda_2(L)} {\tilde{\theta}^{2}}, \frac{2\sqrt{2}(d_{\min})}{\tilde{\theta}}\right).
	\end{equation}
	Let $\alpha\in (0,\alphamaxtilde]$ and let  the step sizes $\bar{\tau},\bar{\beta},\bar{\nu},\bar{\delta}$ be as in Lemma~\ref{lem:stepsizesagg}. Then, for all $\k$, the inclusion in  \eqref{eq:xupdateagg} has a unique solution. Moreover, the sequence $(x^k,s^k,\bs{z}^k,\bs{\lambda}^k)_{k\in\N}$ generated by  Algorithm~\ref{algo:1agg} converges to an equilibrium $(x^*,\1\otimes\avg(x^*)-x^*,\bs{z}^*,\bs{\lambda}^*)$, where  $x^*$ is the \gls{v-GNE} of the game in \eqref{eq:game}. \hfill $\square$
\end{thm}
		\begin{proof1}
			 Similarly to Lemma~\ref{lem:derivation}, we first show that Algorithm~\ref{algo:1agg} can be recast as a \gls{PPPA}, applied to find a zero of the operator $\tilde{\Phi}^{-1}\tilde{	\mc{A}}$. Then, we restrict our analysis to the invariant subspace
		\begin{align}\label{eq:Sigma}
	\Sigma\coloneqq \{(x,s,\bs{v},\bs{\lambda})\in\R^{2n+Em+Nm}\mid \avg(s)=\0_{\bar{n}} \}.
\end{align}
A detailed proof is in Appendix~\ref{app:th:main1agg}. \hfill$\blacksquare$
			\end{proof1}
		\begin{rem1}
		The update in \eqref{eq:xupdateagg} is implicitly defined by a strongly monotone inclusion, or, equivalently, variational inequality (see Appendix~\ref{app:th:main1agg}). We emphasize that there are several iterative methods to find the unique solution (with linear  rate)  \cite[§26]{Bauschke2017} and that,
		as in Remark~\ref{rem:reminexact}, convergence is guaranteed even if the solution is  approximated at each step (with summable errors).
		  \hfill $\square$
	\end{rem1}

	\begin{rem1}
		If, for some $i\in\mc{I}$,
		there exists a function $\varphi_i$ such that $\nabla_{ \! y}\varphi_i(y,s_i^{k+1})=\bs{\tilde{F}}_{\! i}(y,y+s_i^{k+1})$,
		 then the  update of $x_i^k$ in Algorithm~\ref{algo:1agg} can be simplified as
		\begin{align}
		\nonumber	x_i^{k+1}  &= 	\underset{y \in \Omega_i} {\argmin}\bigl(  \varphi_i(y,s_i^{k+1})\textstyle+ \frac{1}{2\alpha\tau_i}  \bigl\|  y-x_i^{k} \bigr\|^2 \\
		& \qquad  \nonumber  \textstyle+\frac{1}{\alpha}(A_i^\top \lambda_i^k)^\top y
		 +\textstyle \frac{1}{\alpha} \textstyle  \bigl( \sum_{j=1}^{N}w_{ij} (\sigma_i^k- \sigma_j^{k}) \bigr)^\top y  \bigr),
		\end{align}
		as in  Lemma~\ref{lem:derivation}. For scalar games (i.e., $\bar{n}=1$)  this condition holds for all $i\in\mc{I}$. Another noteworthy example is  that of a cost  $f_i(x_i,\avg(x))=\bar{f_i}(x_i)+(Q_i\avg(x))^\top x_i$, for some function $\bar{f}_i$ and symmetric matrix $Q_i$, which models  applications as the Nash--Cournot game described in \cite{Koshal_Nedic_Shanbag_2016} and the resource allocation problem considered in \cite{BelgioiosoGrammatico_L_CSS_2017}. In this case,
		$\varphi_i(x_i,s_i)=\bar{f}_i(x_i)+(Q_i (s_i+x_i))^\top x_i-\textstyle \frac{N-1}{2N}x_i^\top Q_i x_i.$
		\hfill$\square$
	\end{rem1}
    \begin{rem1}
    	Inertial/relaxed versions of  Algorithm~\ref{algo:1agg} can be studied as in §\ref{sec:acc}; further, linear convergence can be established for aggregative games without coupling constraints, based on  the restricted strong monotonicity of $\tilde{\Fa}$ (see the proof of Lemma~\ref{lem:amon_agg} in  Appendix~\ref{app:th:main1agg}), as in Theorem~\ref{th:mainNE}. \hfill $\square$
    \end{rem1}
	\section{Numerical simulations}\label{sec:numerics}
	\subsection{Nash--Cournot game}\label{subsec:NC}
	We consider a Nash--Cournot game \cite[§6]{Pavel2018}, where $N$ firms produce a commodity that is sold to $m$ markets. Each firm $i\in\mc{I}=\{1,\dots,N\}$  participates in $n_i\leq m$
	 of the markets, and decides on the quantities $x_i\in\R^{n_i}$ of commodity to be delivered to these $n_i$ markets. The quantity of product that each firm can deliver is bounded by the local constraints $\0_{n_i}\leq x_i\leq X_i$. Moreover, each market $k=1,\dots,m$ has a maximal capacity $r_k$. This results in the shared affine constraint $Ax\leq r$, with $r=\col((r_k)_{k={1,\dots,m}})$ and $A=[A_1 \dots A_N]$, where $A_i\in\R^{m\times n_i}$ is the matrix that expresses which markets firm $i$ participates in. Specifically,  $[A_i]_{k,j}=1$ if $[x_i]_j$ is the amount of product sent to the $k$-th market by agent $i$, $[A_i]_{k,j}=0$ otherwise, for all $j=1,\dots, n_i$, $k={1,\dots,m}$.  Hence, $Ax=\textstyle \sum_{i=1}^N A_ix_i \in \R^m$ is the vector of the quantities of total product delivered to the markets.  Each firm $i$ aims at maximizing its profit, i.e.,  minimizing the cost function
	$J_i(x_i,x_{-i})=10^{-3}*(c_i(x_i)-p(Ax)^\top A_ix_i)$.
	Here, $c_i(x_i)=x_i ^\top Q_i x_i+q_i^\top x_i$ is firm $i$'s production cost, with $Q_i\in \R^{n_i\times n_i}$, $Q_i\succ 0$, $q_i\in \R^{n_i}$. Instead, $p:\R^m\rightarrow \R^m$ associate to each market a price that depends on the amount of product delivered to that market. Specifically, the price for the market $k$, for $k=1,\dots,m$, is $[p(x)]_k=\bar P_k$ -$\chi_k [Ax]_k$, where $\bar P_k$, $\chi_k>0$.

	We set $N=20$, $m=7$. The market structure (i.e., which firms are allowed to participate in which of the $m$ markets) is defined as in \cite[Fig.~1]{Pavel2018}; thus $x=\col((x_i))_{i\in\mc{I}})\in\R^n$ and $n=32$. The firms cannot access the production of all the competitors, but they are allowed to communicate with their neighbors on a randomly generated connected graph. We select randomly with uniform distribution $r_k$ in $[1,2]$, $Q_i$ diagonal with diagonal elements in $[1,8]$, $q_i$ in $[1,2]$, $\bar{P}_k$ in $[10,20]$, $\chi_k$ in $[1,3]$, $X_i$  in $[5,10]$, for all $i\in \mc{I}$, $k=1,\dots,m$.

	The resulting setup satisfies all our theoretical assumptions \cite[§VI]{Pavel2018}.
	We set $\alpha=\alphamax \approx 0.7$ as in  Lemma~\ref{lem:strongmon_constant} and we choose the step sizes as in Lemma~\ref{lem:stepsizes} to satisfy all the conditions of Theorem~\ref{th:main1}.

	We compare the performance of Algorithm~\ref{algo:1} versus that of the pseudo-gradient method  in \cite[Alg.~1]{Pavel2018}, which is to the best of our knowledge the only other available single-layer fixed-step scheme to solve \gls{GNE} problems under partial-decision information. In \cite[Alg.~1]{Pavel2018},  we choose the parameter $c$ that maximize the step sizes $\tau$, $\nu$, $\sigma$, provided that the conditions in \cite[Th.~2]{Pavel2018} are satisfied. This results in very small step sizes, e.g., $\tau^*\approx 10^{-5}.$

\begin{figure}[t]
		\centering
		\includegraphics[width=\columnwidth]{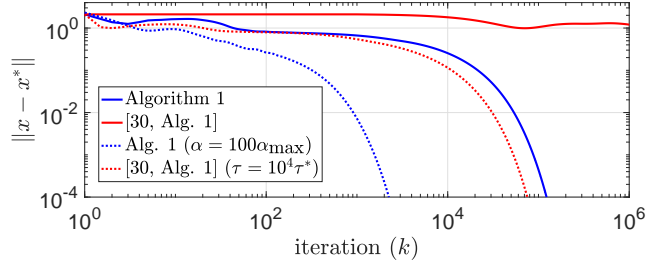}
%
		\caption{Comparison of Algorithm~\ref{algo:1} and \cite[Alg.~1]{Pavel2018} for different parameters (the solid line for the theoretical step sizes).
 }
 \label{fig:1}
	\end{figure}

	The results are illustrated in Figure~\ref{fig:1}, where the two Algorithms are initialized with the same random initial conditions.  \cite[Alg.~1]{Pavel2018} is extremely slow, due to the small step sizes; and our \gls{PPPA} method shows a much faster convergence. According to our numerical experience, the bounds on the parameters are  conservative, and in effect we observe faster convergence for larger step sizes. For  \cite[Alg.~1]{Pavel2018}, the fastest convergence is attained by setting the step sizes $10^4$ times bigger than the theoretical bounds; for larger   steps, convergence is lost.

\begin{figure}[t]
		\centering
		\includegraphics[width=\columnwidth]{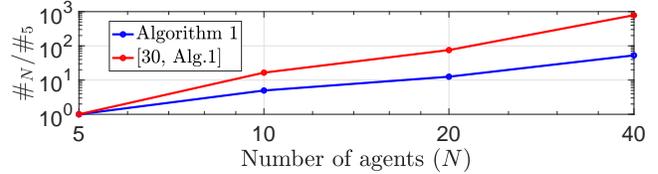}
		\caption{
			Variation of the number of iterations $\#_N$  needed to reach a precision of  $\|x^k-x^*\|\leq 10^{-2}$ for different values of the number of agents $N$ (in logarithmic scale).
		}\label{fig:1b}
	\end{figure}

	We  repeat the simulation for different numbers of agents (and random market structures). Differently from Algorithm~1,  the upper bounds for the step sizes in \cite[Alg.~1]{Pavel2018} decrease when $N$ grows (see §\ref{subsec:rate}), resulting in a greater performance degradation, as shown in Figure~\ref{fig:1b} (with theoretical parameters for our \gls{PPPA}, and steps $10^3$ times larger than their upper bounds for \cite[Alg.~1]{Pavel2018}).

\begin{figure}[t]
		\centering
		\includegraphics[width=\columnwidth]{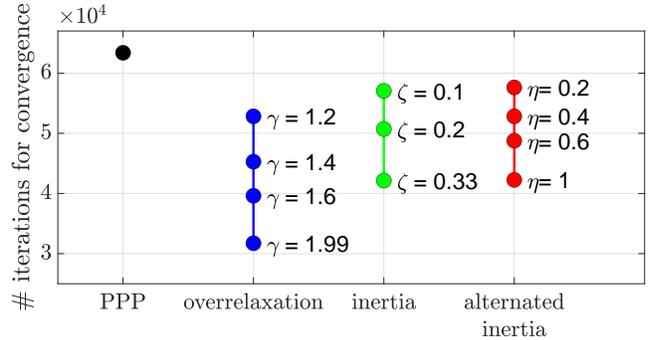}
		\caption{Number of iterations needed to reach a precision of $\|x^k-x^*\|\leq 10^{-2}$, with different acceleration schemes and parameters.}\label{fig:2}
	\end{figure}

	Finally, we apply the acceleration schemes discussed in Section~\ref{sec:acc} to Algorithm~\ref{algo:1},  with parameters that theoretically ensure  convergence.  The impact is remarkable, up to halving the number of iterations needed for  convergence, as shown in Figure~\ref{fig:2}.

	\subsection{Charging of plug-in electric vehicles}\label{subsec:PEV}
	\begin{figure}[t]
		\centering
		\includegraphics[width=\columnwidth]{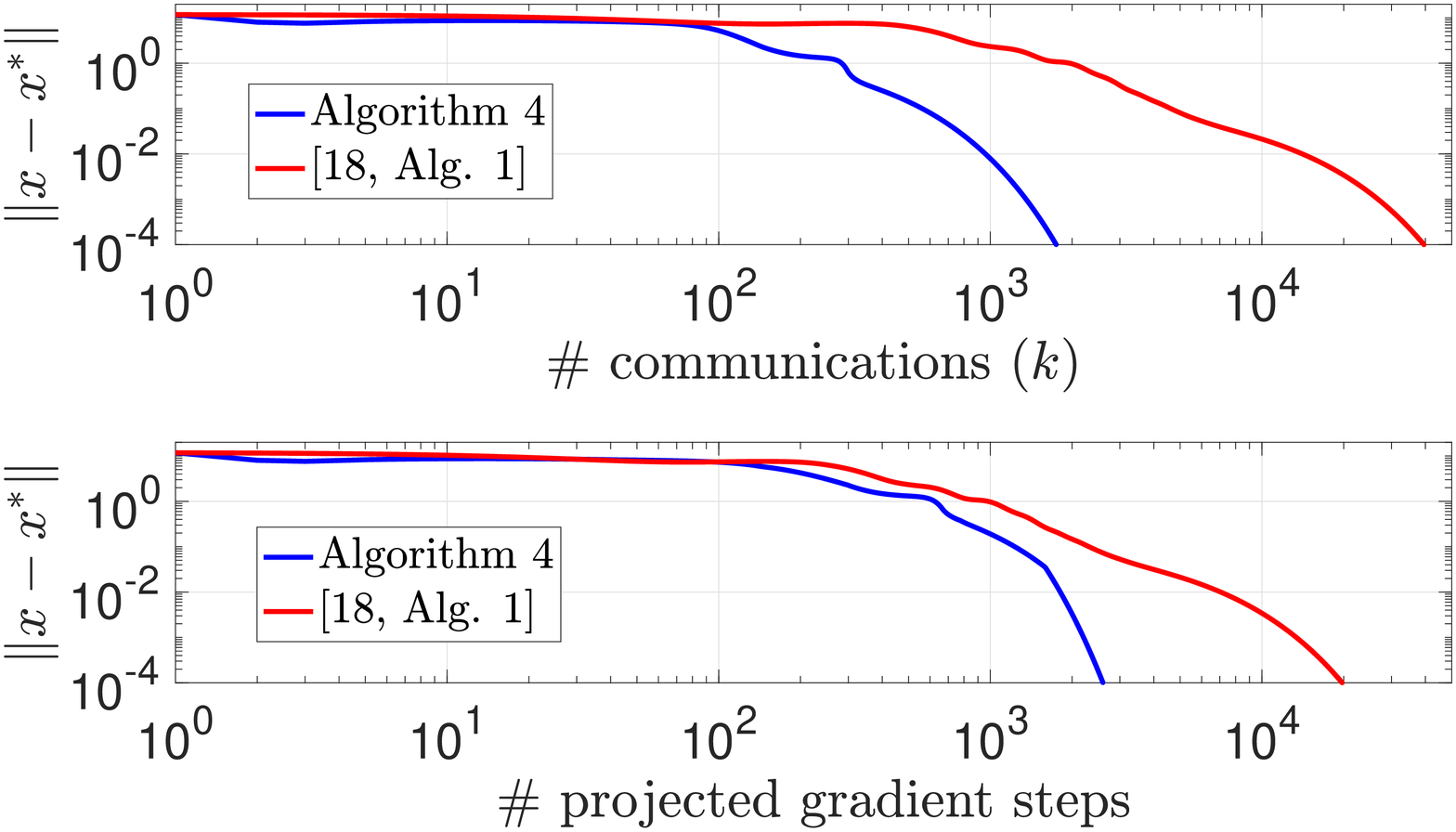}
		\caption{Distance of the primal variable from the \gls{v-GNE}. Algorithm~\ref{algo:1agg} outperforms \cite[Alg.~1]{GadjovPavel_aggregative_2019}, in terms of both  communication rounds and performed projected gradient steps.}\label{fig:3}
		\vspace{1em}
		\includegraphics[width=\columnwidth]{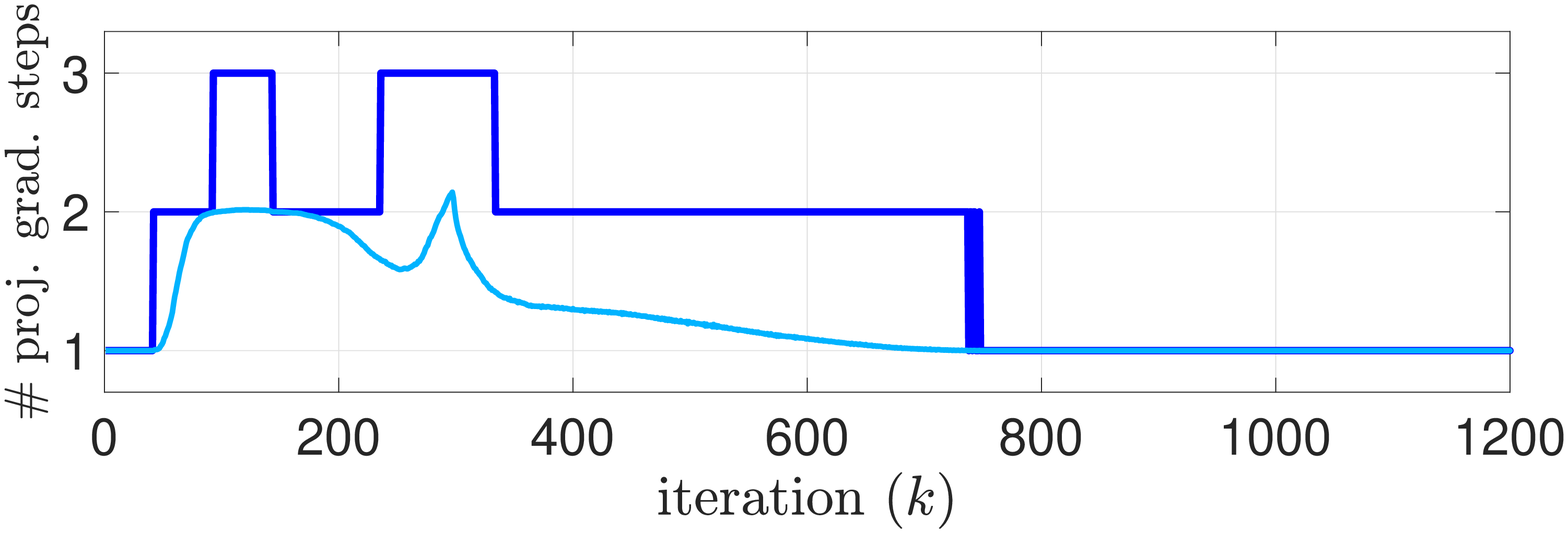}
		\caption{Maximum (blue) and average (light blue) number of projected gradient steps performed by the agents at each iteration in Algorithm~\ref{algo:1agg},  with guaranteed accuracy of  $\varepsilon^k=1/k^2$.}\label{fig:4}
		\vspace{1.5em}
		\includegraphics[width=\columnwidth]{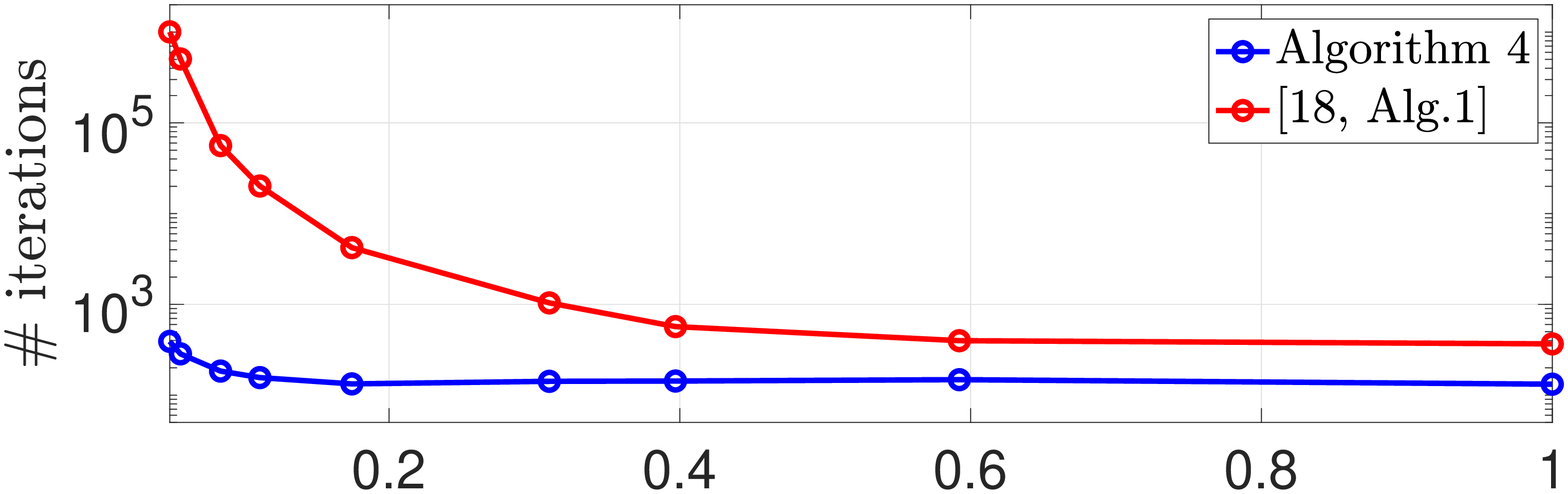}
		\caption{Number of iterations to reach a precision of  $\|x-x^*\|\leq 10^{-2}$ for different values of the algebraic connectivity, where $\uplambda_{2}(L)=1$ indicates a complete graph (all the graphs are doubly stochastic). }\label{fig:5}
	\end{figure}
	We consider the charging scheduling problem for a group of plug-in electric vehicles, modeled by an aggregative game \cite{Grammatico2017}. Each user  $i\in\mc{I}=\{1,\dots,N\}$ plans the charging of its vehicle for an horizon of $24$ hours, discretized into $\bar{n}$ intervals; the goal is to choose the energy injections $x_i\in\R^{\bar{n}}$ of each time interval to minimize its cost $J_i(x_i,\avg(x))=g_i(x_i)+p(\avg(x))^\top  x_i,$
	where  $g_i(x_i)=x_i^\top Q_i x_i +c_i^\top x^i$ is the battery degradation cost, and $p(\xi)=a(\xi+d)+b\1_{\bar{n}}$ is the cost of energy, with $b$  a baseline price, $a$ the inverse of the price elasticity and $d\in\R^{\bar{n}}$
	 the inelastic demand (not related to vehicle charging) along the  horizon. We assume a maximum injection per interval and a desired final charge level for each user, resulting in the local constraints $\Omega_i= \{y\in[0_{\bar{n}},\bar{x}_i]\mid \1_{\bar{n}}^\top y=\gamma_i  \}$. Moreover, we consider the transmission line constraints $\0_{\nbar}\leq \sum_{i\in\mc{I}} x_i \leq \bar{c}N$.

	We set $N=1000$, $\nbar=12$. For all $i\in\mc{I}$, we select with uniform distribution $c_i$ in $[0.55,0.95]$, $Q_i\succ 0$ with diagonal and off-diagonal elements in $[0.2, 0.8]$ and $[0,0.05]$, respectively,  $\gamma_i$ in $[0.6,1]$; $[\bar{x}_i]_j=0.25$ with probability $20\%$, $[\bar{x}_i]_j=0$ otherwise. We set $[\bar{c}]_j$ as $0.04$ if $j\in\{1,2,3, 11,12 \}$, as $0.01$ otherwise (corresponding to more restrictive limitations in the daytime); $a=0.38$, $b=0.6$ and $d$ as in  \cite{Grammatico2017}.
	We check numerically that Standing Assumptions~\ref{Ass:Convexity},~\ref{Ass:StrMon} hold, and  let the agents communicate over a randomly generated  connected graph.
	We implement Algorithm~\ref{algo:1agg}, by performing only a finite number of gradient steps per iteration; each agent uses the stopping criterion in \eqref{eq:stopping criterion} to ensure an accuracy of $\varepsilon^k=1/k^2$.
	Figure~\ref{fig:3} compares the performance of Algorithm~\ref{algo:1agg}
	and \cite[Alg.~1]{GadjovPavel_aggregative_2019} {(which requires two rounds of communication per iteration)}, with  step sizes set to their theoretical upper bounds. Notably, our \gls{PPPA} significantly outperforms \cite[Alg.~1]{GadjovPavel_aggregative_2019}, even in terms of total projected gradient steps required  (for Algorithm~\ref{algo:1agg}, we consider the maximum per iteration). Interestingly, Figure~\ref{fig:4} shows that the maximum number of performed gradient steps at each iteration is $3$ and  decreases as the iteration converges, despite the increasing accuracy required in the local optimizations (see also Remark~\ref{rem:reminexact}).

		Differently from our \gls{PPPA},  the upper bounds for the step sizes in  \cite{GadjovPavel_aggregative_2019} are  proportional to the quantity $\mu_{\tilde{A}}$ in \cite[Lem.~4]{GadjovPavel_aggregative_2019}, hence they depend  on  $\uplambda_2(L)$, $\theta_0$, $\mu$, $\theta$ (but not on $N$, cf. §~\ref{subsec:NC}, \ref{subsec:rate}); in turn, we expect these parameters to affect to a larger extent the convergence speed for the \gls{FB} method. In Figure~\ref{fig:5} we compare the two algorithms, with $N=10$,  for different values of the communication graph connectivity: in the considered range, the number of iterations to converge varies by a factor $2$ for Algorithm~\ref{algo:1agg}, by a factor $10^3$ for \cite[Alg.~1]{GadjovPavel_aggregative_2019}.
	~
	\section{Conclusion}
	Inexact preconditioned proximal-point methods are extremely efficient to design fully-distributed single-layer generalized Nash equilibrium seeking algorithms. The advantage is that  convergence  can be guaranteed for much larger step sizes compared to pseudo-gradient-based algorithms.  In fact, in our numerical experience, our algorithms proved much faster than the existing methods, resulting in a considerable reduction of communication and computation requirements. Besides, our operator-theoretic approach facilitates the design of acceleration schemes, also in the partial-decision information setup.
	As future work, it would be highly  valuable to relax our monotonicity and connectivity assumptions, namely to allow for merely monotone game mappings  and  jointly connected networks, and to address the case of nonlinear coupling constraints.
	\appendix

\section{Proof of Theorem~\ref{th:FNEconvergence}}\label{app:th:FNEconvergence}
For all $\k$, let $z^{k}\coloneqq \omega^k+\gamma^k(u^k-\omega^k)$, so that $\omega^{k+1}=z^{k}+\gamma^k e^k$. Consider any  $\o^*\in C$.  We have, for all $\k$,
	\begin{align}
   & \hphantom{{}={} } \| z^{k}-\omega^*\|^2_P
	\nonumber
	\\ \nonumber & = \|\omega^k-\omega^*\|_P^2-2\gamma^k\langle \omega^k-u^k \mid \omega^k-\omega^*\rangle_P
	 \\	\nonumber  & \hphantom{ {}={} \|\omega^k-\omega^*\|_P^2}
	  + (\gamma^k)^2 \|u^k-\omega^k\|^2_P  \\
	\label{eq:proof1}	& \leq \|\omega^k-\omega^*\|_P^2-\gamma^k(2-\gamma^k)\|u^k-\omega^k\|_P^2,
	\end{align}
	where the inequality follows by Lemma~\ref{lem:FQNE}.
	\newline
	{\it (i)} By \eqref{eq:proof1}, $\|z^{k}-\omega^*\|_P \leq \|\omega^k-\omega^*\|_P$, and the conclusion follows by the Cauchy--Schwartz inequality. \newline
	{\it (ii)} By $(\gamma^k\|e^k\|_P)_{k\in\N}\in\ell^1$ and point {\it (i)}, $(\omega^k)_{k\in\N}$ is bounded. Let $c\coloneqq \sup_{\k}\|\omega^k-\omega^*\|_P<\infty$ and   $\epsilon^k\coloneqq 2c(\gamma^k \|e^k\|_P)+(\gamma^k \|e^k\|_P)^2$, for all $\k$. Clearly, $(\epsilon^k)_{\k}\in\ell^1$. Moreover, for all $\k$ we have
	\begin{align}
	\nonumber & \hphantom{{}\leq{}}\|w^{k+1}-w^*\|_P^2  \\
	\nonumber &\leq ( \|z^{k}-w^*\|_P+\gamma^k \|e^k\|_P)^2 \\
	\label{eq:proof2}  & \leq \|w^k-w^*\|_P^2-\gamma^k(2-\gamma^k)\|u^k-w^k\|_P^2+\epsilon^k,
	\end{align}
	and the thesis follows by recursion. \newline
	{\it (iii)} By  \eqref{eq:proof2}, \cite[Prop.~3.2(i)]{Combettes2001} and \cite[Th.~3.8]{Combettes2001}. \newline
	{\it (iv)} By definition of resolvent, $\omega^k-u^k\in\mc{B}(u^k)$; hence
	 \begin{equation}
	 \label{eq:proof3}
		\langle  u^k-\omega^* \mid \omega^k-u^k \rangle_P \geq \mu_{\mc{B}}\|u^k-\omega^*\|_P^2. \end{equation}
		By the Cauchy--Schwartz inequality,   $\|\omega^k-u^k\|_P\geq \mu_{\mc{B}}\|u^k-\omega^*\|_P$. Thus, \eqref{eq:proof3} yields
	\begin{align}
	\nonumber & \hphantom{{} \leq {}}	\| \omega^k-\omega^*\|_P^2 \\
	\nonumber & = \|u^k-\omega^*\|_P^2+2\langle  u^k-\omega^* \mid \omega^k-u^k \rangle_P +\| \omega^k-u^k\|_P^2 \\
	\label{eq:proof4} &\geq  (1+\mu_{\mc{B}})^2  \|u^k-\omega^*\|_P^2.
	\end{align}
	If $\gamma^k \leq 1$,  by the Cauchy--Schwartz inequality and \eqref{eq:proof4}, we have  $\|z^{k}-\omega^*\|_P \leq (1-\gamma^k)\| \omega^k-\omega^*\|_P+\gamma^k\|u^k-\omega^*\|_P\leq (1-\frac{\gamma^k \mu_{\mc{B}}}{1+\mu_{\mc{B}}})\| \omega^k-\omega^*\|_P$. For $\gamma^k>1$, we can write
	\begin{align}
	& \hphantom{{}={}} \| z^{k}-\omega^*\|_P^2  \nonumber
	\\ \nonumber
	& = (1-\gamma^k)^2\| \omega^k-\omega^*\|_P^2 +\gamma^k(2-\gamma^k) \|u^k-\omega^*\|_P^2
	\\ \nonumber
	&
	 \hphantom{{}={}} +2\gamma^k(1-\gamma^k) \langle u^k-\omega^* \mid \omega^k-u^k \rangle_P
	\\ \nonumber
	& \leq  (1-\gamma^k)^2\| \omega^k-\omega^*\|_P^2
	\\
	&  \hphantom{{}={}} +\gamma^k (2(1+\mu_{\mc{B}} ) -\gamma^k(1+2\mu_{\mc{B}})) \|u^k-\omega^*\|_P^2 \label{eq:proof5} \\
	& \leq (\max(1- \textstyle \frac{\gamma^k \mu_{\mc{B}}}{1+\mu_{\mc{B}}},\gamma^k-1))^2\|\omega^k-\omega^*\|_P^2, \label{eq:proof6}
	\end{align}
	where the first equality follows by rearranging the terms in \eqref{eq:proof1}; in the first inequality we used \eqref{eq:proof3}; the last inequality  follows by taking into account that the second term in \eqref{eq:proof5} is nonpositive if $\gamma^k \in (1,1+\textstyle\frac{1}{1+2\mu_{\mc{B}}}]$ and can be upper bounded via \eqref{eq:proof4} if $\gamma^k \in [1+\textstyle\frac{1}{1+2\mu_{\mc{B}}},2)$.
	 Finally, assume that $\o^k\in C$, and choose $u^k=\omega^k$. Then \eqref{eq:proof6} implies $\o^k= \o^*$, hence $C$ must be a singleton. \hfill $\blacksquare$

	\section{Proof of Theorem~\ref{th:main1agg}}
	\label{app:th:main1agg}
	Analogously to Lemma~\ref{lem:derivation}, it can be shown that Algorithm~\ref{algo:1agg} is equivalent to the iteration
	\begin{equation}\label{eq:algcompactagg}
	\obs^{k+1}\in \J_{\tilde\Phi^{-1}\tilde\Ac}(\obs^k), \quad  \obs^0=\bar{\obs}^0,
	\end{equation}
	where $\bar{\obs^0}=(\xbs^0,\0_n,\0_{Em},\lbs^0)$, for some $\xbs^0\in\Omega$, $\lbs^0\in \R_{\geq 0}^{Nm}$, modulo the transformation $\bs{z^k}=\bs{V}_{\!\!m}^\top \bs{v^k}$.

	First, we show that the iteration in \eqref{eq:algcompactagg} is uniquely defined. For all $i\in\mc{I}$, let $\Fc_i(y,\vartheta^k)\coloneqq \alpha\bs{\tilde{F}}_{\! i}(y,y+s_i^{k+1})
	+\frac{1}{\tau_i}(y- x_i^k)+A_i^\top {\lambda}_i^k +\textstyle \sum_{j=1}^{N}w_{i,j}(\sigma_i^k-\sigma_j^k)+\nc_{\Omega_i}(y)$, where $\vartheta^k=(x^k,s^{k+1},s^k,\bs{\lambda}^k)$. We note that ${\bs{\tilde{F}}}_{\! i}$ is $\tilde{\theta}$-Lipschitz, because $\bs{\tilde{F}}$ is $\tilde{\theta}$-Lipschitz by Lemma~\ref{lem:LipschitzExtPseudoagg}. Then, by monotonicity of the normal cone, we have  $\langle y-y^\prime \mid \Fc_i(y,\vartheta^k)-\Fc_i(y^\prime,\vartheta^k)\rangle \geq (\tau_i^{-1}-\alpha\sqrt{2}\tilde\theta)\|y-y^\prime\|^2$, for any $y,y^\prime\in\R^n$, for any $\vartheta^k$. By the assumption on $\alpha$, $\Fc_i$ is strongly monotone in $y$ for any $\vartheta^k$, hence the inclusion in  \eqref{eq:xupdateagg} has a unique solution, for any $\vartheta^k$ \cite[Cor.~23.37]{Bauschke2017}.
	Therefore, it also holds that $\dom(\J_{\tilde\Phi^{-1}\tilde\Ac})=\R^{2n+Em+Nn}$ and that
	$\J_{\tilde\Phi^{-1}\tilde\Ac}$ is single-valued.

	We  turn our attention to the set $\Sigma$ in \eqref{eq:Sigma}.  As in Remark~\ref{rem:invariance},  for any $\varsigma \in \Sigma$, $J_{\tilde\Phi^{-1}\tilde\Ac}(\varsigma)\in \Sigma$; hence $\Sigma$ is invariant for \eqref{eq:algcompactagg}.
    Moreover, $\obs^{0}\in\Sigma$. Hence, in \eqref{eq:algcompactagg}, it is enough to consider  the operator $J_{\tilde\Phi^{-1}\tilde\Ac}\!\!\mid_\Sigma:\Sigma\rightarrow \Sigma$, where $\mc{B}|_{\Sigma}$ is the restriction of the operator $\mc{B}$ to  $\Sigma$, i.e., $\mc{B}|_\Sigma(\obs)=\mc{B}(\obs)$ if $\obs\in \Sigma$, $\mc{B}|_\Sigma(\obs)=\varnothing$ otherwise. By invariance and \eqref{eq:explicitinclusion}, it also follows that $J_{\tilde\Phi^{-1}\tilde\Ac}\!\mid_\Sigma=J_{\tilde\Phi^{-1}\tilde\Ac\mid_\Sigma} \! \mid_\Sigma $. Thus, the iteration in \eqref{eq:algcompactagg} is rewritten as
    \begin{equation}\label{eq:algcompactagg2}
    \obs^{k+1}= \J_{\tilde\Phi^{-1}\tilde\Ac\mid_\Sigma}(\obs^k), \quad  \obs^0=\bar{\obs}^0.
    \end{equation}
We show the convergence of \eqref{eq:algcompactagg2} by studying   the properties of $\tilde{\Ac}\!\mid_\Sigma$.
We start by characterizing the zero set.

	\begin{lem1}\label{lem:zerosagg}
		The following statements holds:
		\begin{itemize}[topsep=-2em]
			\item [(i)] If $\col(x^*,s^*,\bs{v}^*,\bs{\lambda}^*)\in\zer(\tilde{\Ac}\!\mid_\Sigma)$, then $s^*=\1_N\otimes\avg(x^*)-x^*$  and $x^*$ is the \gls{v-GNE} of the game in \eqref{eq:game}.
			\item[(ii)] 	$\zer(\tilde{\Ac}\!\mid_\Sigma) \neq \varnothing$.  \hfill $\square$
		\end{itemize}
	\end{lem1}

	\begin{proof1} Let  $\bs{V}_{\!\!q}\coloneqq V\otimes I_q$, $\bs{L}_q\coloneqq L\otimes I_q=\bs{V}_{\!\!q}^\top\bs{V}_{\!\!q}$, for any $q>0$;  hence, under Standing Assumption~\ref{Ass:Graph}, we have
			\begin{alignat}{3}
			\label{eq:null_Lq}
			\Null \left(\bs{L}_q\right)   & =\Null \left(\bs{V}_{\!\!q}\right)     &   & = \Range(\1_N\otimes I_q)
			\\
			\label{eq:range_Lq}
			\!\!\!\!\!\! \Range(\bs{V}^\top_{\!\!q}  ) & \supseteq \Range\left(\bs{L}_q \right) &   & =\Null(\1_N^\top\otimes I_q).
			\end{alignat}
			\textit{(i)} Let us consider any $ {\obs^*}=\operatorname{col}({ { {x}}^*},s^*,  {\bs v}^*,{\bs{ \lambda}}^*) \in{\zer(\tilde{\Ac}\!\mid_\Sigma)}$, and let  $\sigma^*=x^*+s^*$; then we have
			\begin{subequations}
				\label{eq:eqcond}
				\begin{align}
				\label{eq:eqconda}
				\0_{\bar{n}} & \in \alpha\bs{\tilde{F}}( {x}^*, {\sigma}^*)+\bs{L}_{\bar{n}} {\sigma^*}+\nc_{{\Omega}}( { {x}^*}) +
				\bs{A}^\top {\bs{\lambda}}^*
				\\
				\label{eq:eqcondb}
				\0_{\bar{n}} & = \bs{L}_{\bar{n}}{{\sigma}^*}
				\\
				\label{eq:eqcondc}
				\0_{Em}      & = -\bs{V}_{\!\!m} {\bs{\lambda}^*}
				\\
				\label{eq:eqcondd}
				\0_{Nm}      & \in \bs b+ \nc_{\R^{Nm}_{\geq 0}}({\bs{\lambda}}^*)  {-\bs{A}{x}^* + {\bs{V}_{\!\!m}^\top {\bs{v}}^*
				}}
				\end{align}
			\end{subequations} By \eqref{eq:eqcondc} and by \eqref{eq:null_Lq}, we have $ {\bs{\lambda}}^*=\1_N\otimes \lambda^*$, for some $\lambda^*\in \R^m$; by \eqref{eq:eqcondb} and since ${\obs}^*\in \Sigma$, it must hold $\sigma^*=x^*+s^*=1_N\otimes \avg(x^*)$. It is then enough to prove that the pair $(x^*, \lambda^*)$ satisfies the \gls{KKT} conditions in \eqref{eq:KKT}.
			By \eqref{eq:eqconda},  by recalling that $\bs{A}^\top (\1_N\otimes \lambda^*)=A^\top \lambda^*$ and
			$\bs{\tilde{F}}(x^*,\1_N \otimes x^*)=F(x^*)$, we retrieve the first KKT condition in \eqref{eq:KKT}.
			We obtain the second \gls{KKT} condition by  left-multiplying both sides of \eqref{eq:eqcondd} with
			$(\1_N^\top \otimes I_m)$
			and using that
			$(\1_N^\top \otimes I_m)\bs{b}=b$,
			$(\1_N^\top \otimes I_m)\bs{L}_{m}=0$
			by \eqref{eq:null_Lq} and symmetry of $L$,
			$\textstyle{(\1_N^\top \otimes I_m)\bs{A}=A}$ and
			$\textstyle {(\1_N^\top \otimes I_m)N_{\R_{\geq 0}^{Nm}}(\1_N\otimes\lambda^*)=N\mathrm{N}_{\R_{\geq 0}^{ m}}(\lambda^*)=\mathrm{N}_{\R_{\geq 0}^{ m}}(\lambda^*)}$.
			\textit{(ii)}  	Let us consider any pair  $(x^*,\lambda^*)$ satisfying the \gls{KKT} conditions in \eqref{eq:KKT} (one such pair exists by  Assumption~\ref{Ass:StrMon}). We next show that there exists $\bs{v}^*\in \R^{Em}$ such that  $\obs^*=\col(x^*,\1_N\otimes \avg(x^*)-x^*,\bs{z}^*,\1_N\otimes\lambda^*)\in \zer({\tilde{\Ac}\!  \mid_\Sigma})$.  Clearly, $\obs^*\in \Sigma$. Besides, $\obs^*$  satisfies the conditions \eqref{eq:eqconda}-\eqref{eq:eqcondc}, as in point {\it(i)}.
			By \eqref{eq:KKT},  there exists $u^*\in \mathrm{N}_{\R^m_{\geq 0}}(\lambda^*)$ such that  $Ax^*-b-u^*=\bs 0_n$.
			Also,  $\mathrm{N}_{\R_{\geq 0}^{N m}}(1_N\otimes\lambda^*)=\prod_{i\in \mc I}\mathrm{N}_{\R^m_{\geq 0}}(\lambda^*)$, and it follows by properties of cones that $\operatorname{col}\left(u_1^*,\dots,u_N^*\right)\in \mathrm{N}_{\R_{\geq 0}^{N m}}(\1_N\otimes\lambda^*)$, with $ u_1^*= \dots =u_N^*=\frac{1}{N}u^*$. Hence
			$(\1_{N}^\top \otimes I_{m})\left( -\bs{A} x^*+\bs{b}  +\operatorname{col}\left(u_1^*,\dots,u_N^*\right)\right)=b- Ax^*+u^*=\bs 0_m,$ or $  -\bs{A} x^*+\bs{b}  +\operatorname{col}\left(u_1^*,\dots,u_N^*\right)\in \Null (\1_{N}^\top \otimes I_{m})\subseteq \Range (\bs{V}^\top_{\!\!m})$, by \eqref{eq:range_Lq}. Therefore there exists $\bs{v}^*$ such that also the condition \eqref{eq:eqcondd} is satisfied, for which $\obs^*\in\zer(\tilde{\Ac})$.
			\hfill  $\blacksquare$
	\end{proof1}

	Next, similar to Lemma~\ref{lem:strongmon_constant}, we show  restricted monotonicity of the operator $\tilde{\Ac}\!\mid_\Sigma$.

	\begin{lem1}\label{lem:amon_agg}
		Let $\alpha\in(0,{\tilde{\alpha}_{\textnormal{max}}}]$, with ${\tilde{\alpha}_{\textnormal{max}}}$ as in \eqref{eq:alphamaxagg}. Then $\tilde{\Ac}\!\mid_\Sigma$ is restricted monotone.
	\end{lem1}

	\begin{proof1} The operator $\tilde{\Ac}\!\mid_\Sigma$  is the sum of three components, as in \eqref{eq:opAtilde}. The third is monotone by properties of the normal cones \cite[Th.~20.25]{Bauschke2017}, the second because it is a linear skew-symmetric operator \cite[Ex.~20.35]{Bauschke2017} (and restriction does not cause loss of monotonicity, by definition).
	For the first term, let $(\obs,\ubs)\in\gra(\tilde{\Ac}\!\mid_\Sigma)$,
	$\obs\coloneqq \col(x,s,\bs{v},\bs{\lambda})$, $\obs^*=\col(x^*,s^*,\bs{v}^*,\bs{\lambda}^*)\in \zer(\tilde{\Ac} \! \mid_\Sigma)$, $\sigma=x+s$, $\sigma^*=s^*+x^*$.
		By Lemma~\ref{lem:zerosagg}, $s^*=\1_N\otimes \avg(x^*)-x^*$. Then, by \cite[Lemma~4]{GadjovPavel_aggregative_2019}, there is a $\tilde{\mu}>0$ such that
				$
        \langle \col(x-x^*,s-s^*) \mid \tilde{\Fa}(x,s) -\tilde{\Fa}(x^*,s^*)\rangle
		 =\langle x-x^* \mid \alpha\bs{\tilde{F}}(x,{\sigma})-	\alpha\bs{\tilde{F}}(x^*,{{\sigma}^*})\rangle
		+{ \langle {\sigma}-{{\sigma}^*} \mid \bs{L}_{\bar{n}}({\sigma}-{{\sigma}^*)}\rangle }
		\geq \tilde{\mu} \| \col(x-x^*,\sigma-\sigma^*)\|^2 \geq \mu_{\tilde{\Fa}}  \|\col(x-x^*,s-s^*)\|^2$, where $\mu_{\tilde{\Fa}}\coloneqq (3-\sqrt{5})\tilde{\mu}/2$ and the last inequality follows by definition of $\sigma$ and bounds on quadratic forms.\hfill $\blacksquare$
	\end{proof1}

	Finally, the preconditioning matrix $\tilde \Phi$ is positive definite by Lemma~\ref{lem:stepsizesagg}. As in Lemma~\ref{lem:strmoninPhi}, by Lemma~\ref{lem:amon_agg}, it holds that $\tilde{\Phi}^{-1}\tilde{\Ac}\!\mid_\Sigma$ is restricted monotone in $\H_{\tilde \Phi}$. In view of \eqref{eq:algcompactagg2}, the conclusion follows analogously to Theorem~\ref{th:main1}. \hfill $\blacksquare$
	\bibliographystyle{plain}
\bibliography{library}

\end{document}